\documentclass[12pt]{amsart}

\usepackage[latin1]{inputenc}
\usepackage{mathptmx}
\usepackage{amscd}

\theoremstyle{theorem}
\newtheorem{theorem}{Theorem}[section]
\newtheorem{lemma}[theorem]{Lemma}
\newtheorem{corollary}[theorem]{Corollary}
\newtheorem{proposition}[theorem]{Proposition}

\theoremstyle{definition}
\newtheorem{remark}[theorem]{Remark}
\newtheorem{definition}[theorem]{Definition}
\newtheorem{example}[theorem]{Example}

\theoremstyle{plain}
\newtheorem*{thm*}{Theorem}

\newcommand{\mm}{\mathfrak m}

\newcommand{\pp}{\mathfrak p}
\newcommand{\qq}{\mathfrak q}
\def\NZQ{\Bbb}

\def\ZZ{{\NZQ Z}}
\def\RR{{\NZQ R}}

\def\NN{{\NZQ N}}
\newcommand{\Mcc}{\mathcal{M}}
\newcommand{\Ncc}{\mathcal{N}}
\newcommand{\Wcc}{\mathcal{W}}

\let\Bbb=\mathbb
\let\phi=\varphi

\def\relint{\operatorname{relint}}

\def\nat{\operatorname{nat}}

\def\Hom{\operatorname{Hom}}
\def\Ext{\operatorname{Ext}}
\def\rank{\operatorname{rank}}

\def\Ass{\operatorname{Ass}}

\def\Ker{\operatorname{Ker}}

\def\lcm{\operatorname{lcm}}

\def\Coker{\operatorname{Coker}}

\def\Im{\operatorname{Im}}

\def\cha{\operatorname{char}}

\let\oldbigwedge\bigwedge
\def\BIGwedge{{\textstyle\oldbigwedge}}
\def\medwedge{{\scriptstyle\oldbigwedge}}
\def\bigwedge{\mathchoice{\BIGwedge}{\BIGwedge}{\medwedge}{}}

\let\iso=\cong

\let\epsilon=\varepsilon

\DeclareMathOperator{\cn}{cn}

\DeclareMathOperator{\Supp}{Supp}
\DeclareMathOperator{\fan}{Fan}
\DeclareMathOperator{\str}{star}

\DeclareMathOperator{\lpnt}{\hbox{\large\bf.}}
\DeclareMathOperator{\pnt}{\raise 0.5mm \hbox{\large\bf.}}

\title{On toric face rings}

\author{Bogdan Ichim}
\address{Universität Osnabrück, FB Mathematik/Informatik, 49069
Osnabrück, Germany;
\phantom{iii} Institute of Mathematics, C.P. 1-764, 70700 Bucharest, Romania}
\email{bogdan.ichim@math.uos.de, bogdan.ichim@imar.ro}

\author{Tim Römer}
\address{Universität Osnabrück, FB Mathematik/Informatik, 49069 Osnabrück, Germany}
\email{troemer@uos.de}

\begin{document}
\begin{abstract}
Following a construction of Stanley we consider
toric face rings associated to rational pointed fans.
This class of rings is a common generalization of the
concepts of Stanley--Reisner and affine monoid algebras.
The main goal of this article is to unify parts of the theories of Stanley--Reisner- and affine monoid algebras.
We consider (non-pure) shellable fan's and the Cohen--Macaulay property.
Moreover, we study
the local cohomology, the canonical module
and the Gorenstein property of a toric face ring.
\end{abstract}
\maketitle

%
%
%
\section{Introduction}
\label{intro}

In this article $K$ always denotes a field.
Several types of $K$-algebras appear naturally in combinatorial commutative algebra.
On the one hand let $\Delta$ be a simplicial complex on the vertex
set $[n]=\{1,\dots,n\}$, i.e.\ $\Delta$ is a set of subsets of $[n]$
such that for $F \subseteq G \in \Delta$ we have $F \in \Delta$.
Consider the Stanley--Reisner ring $K[\Delta]$ associated
to $\Delta$.
The relationship between simplicial complexes and Stanley--Reisner rings
has important applications especially in algebraic combinatorics.
Since the first papers by Hochster, Reisner and Stanley
(see e.g.\ \cite{HO77}, \cite{RE76} and \cite{ST75}),
many people studied this subject. For a detailed discussion
see Bruns--Herzog \cite{BH} or Stanley \cite{ST05}.
On the other hand let $M$ be an affine monoid, i.e.\
$M$ is a finitely generated commutative monoid
which can be embedded into $\ZZ^d$ for some $d \in \NN$.
Let $K[M]$ be the affine monoid ring associated to $M$.
Also the relationship between algebraic properties of $K[M]$
and monoid properties of $M$ was intensively studied in many papers,
which is e.g.\ of particular interest for the theory of toric
varieties. See Bruns--Herzog \cite{BH} or
Bruns--Gubeladze \cite{BRGUBbook}
for details and known results.

Usually Stanley--Reisner rings and affine monoid algebras
are studied separately and there is no deep connections between
these two theories. The main goal of this paper is to present a class
of $K$-algebras which include Stanley--Reisner rings and affine monoid algebras,
and to study algebraic properties of these  rings motivated by known
results from the two special cases.
For this we recall a construction of Stanley in \cite{ST87}.
Let $\Sigma$ be a rational pointed fan in $\RR^d$,
i.e.\ $\Sigma$ is a finite collection of rational pointed cones in $\RR^d$
such that for $C' \subseteq C$ with $C\in \Sigma$
we have that $C'$ is a face of $C$ if and only if $C' \in \Sigma$,
and if $C,C' \in \Sigma$, then $C\cap C'$ is a common face of $C$ and $C'$.
The {\em toric face ring $K[\Sigma]$} of $\Sigma$ over $K$ is defined as follows.
As a $K$-vectors space $K[\Sigma]$ has one basis element $x^a$ for
every $a \in \bigcup_{C\in \Sigma} C \cap \ZZ^d$.
Multiplication in $K[\Sigma]$ is defined by the following rule:
$$
x^a \cdot x^b =
\begin{cases}
x^{a+b}, & \text{if $a$ and $b$ are elements of a common face $C\in \Sigma$;}\\
0, & \text{otherwise.}
\end{cases}
$$
Then $K[\Sigma]$ is naturally a $\ZZ^d$-graded $K$-algebra.
There are several other descriptions of these type of rings.
E.g.\ consider the affine monoid rings $K[C \cap \ZZ^d]$
and the natural face projections $K[C \cap \ZZ^d] \to K[C' \cap \ZZ^d]$
for $C,C' \in \Sigma$ such that $C'$ is a face of $C$.
Then one can show that $K[\Sigma]$ is the inverse limit
$\underleftarrow{\lim} K[C\cap \ZZ^d]$ over that systems of rings.
This was the main point of view in \cite{BBR05}
where the local cohomology groups of rings of such type were studied systematically.
A presentation of the toric face ring $K[\Sigma]$ was computed besides other things in \cite{BRRO05},
and initial ideals of the presentation ideals were considered in \cite{BRRO04}.

There exist two extremal cases.
If the fan is the face poset $\fan(C)$ of a single rational pointed cone $C$,
then $K[\fan(C)]=K[C \cap \ZZ^d]$ is just a normal affine monoid ring.
On the other hand, if $\Sigma$ is a {\em simplicial} fan, i.e.\
all cones $C \in \Sigma$ are simplicial, then
$\Sigma$ is as a poset isomorphic to an abstract simplicial complex with
Stanley--Reisner ring $K[\Sigma]$.
In fact it is easy to see that
every Stanley--Reisner ring appears  as a toric face ring
as was observed in \cite[Example 3.2]{BBR05}.
Thus these example show that toric face rings are indeed a common generalization of the
concepts of Stanley--Reisner and affine monoid rings.
Further examples of toric face rings are
polytopal algebras associated to (embedded) polytopal complexes as studied
by Bruns--Gubeladze in \cite{BRGU01} or \cite{BRGUBbook}.
In the case that the fan is the subfan of $\fan(C)$
for one cone $C$
toric face rings are studied by Yanagawa in \cite{YA06}.

This paper is organized as follows.
In Section \ref{toricfacering} we construct a slightly more general
toric face ring $K[\Mcc_\Sigma]$
for a {\em monoidal complex}
$\Mcc_\Sigma$,
which is a set of affine monoids $\{M_C \subseteq \ZZ^d : C \in \Sigma \}$
such that
$\cn(M_C)=C$ for $C \in \Sigma$,
and if $C,C' \in \Sigma$, $C' \subseteq C$, then $M_{C'} = M_C \cap C'$.
Choosing the monoids $N_C=C \cap \ZZ^d$ for $C \in \Sigma$
gives rise to the toric face ring $K[\Sigma]$
of $\Sigma$. We also determine
the $\ZZ^d$-graded prime and radical ideals of this class of toric face rings.
Many of our results in this paper hold for these type of $K$-algebras.
Especially one has not to require that the monoid rings are normal.
Observe that is possible to extend our constructions in various way
and some of our results hold more generally.
One could consider non-embedded monoids glued together corresponding
to the given fan, or to consider other structures like a poset
instead of a fan. Since we searched for a class of rings which include the classical ones,
and we get nice and clean ``embedded results'' using fans, we restricted ourself to this situation.

Analogously to the notion of (non-pure) shellable simplicial complexes,
constructions of Björner--Wachs \cite{BJWA96, BJWA97}
lead to the definition of (non-pure) shellable fans. See Section \ref{shellfan} for details.
Theorems \ref{shellcm} and \ref{shellseqcm}
in Section \ref{shellfan} show that
well-known results for Stanley-Reisner rings generalize to
toric face rings. More precisely, we can prove:

\begin{thm*}
Let $\Mcc_\Sigma$ be a monoidal complex supported by a
rational pointed fan $\Sigma$ in $\RR^d$ such that
$K[M_C]$ is Cohen--Macaulay for all $C \in \RR^d$.
\begin{enumerate}
\item
If $\Sigma$ is shellable, then
$K[\Mcc_\Sigma]$ is Cohen--Macaulay.
\item
If $\Sigma$ is non-pure shellable, then
$K[\Mcc_\Sigma]$ is sequentially Cohen--Macaulay.
\end{enumerate}
In particular, (i) and (ii) can be applied to the toric face ring $K[\Sigma]$ of $\Sigma$.
\end{thm*}

In Section \ref{localcoho} we consider the local cohomology groups $H^i_{\mm}(K[\Mcc_\Sigma])$
of
$K[\Mcc_\Sigma]$ with respect to its unique $\ZZ^d$-graded maximal ideal $\mm=(x^a : 0\neq a \in |\Mcc_\Sigma|)$,
where $|\Mcc_\Sigma|=\bigcup_{C\in \Sigma} M_C$.
We present explicit complexes which compute the local cohomology groups $H^i_{\mm}(K[\Mcc_\Sigma])$.
In particular, in Corollary \ref{hochster}
we show that
one can obtain Hochster type formulas for the local cohomology of the toric face ring
$K[\Sigma]$,
and we get a direct proof of a variation of the main result of \cite{BBR05} applied to toric face rings
avoiding the technical machinery developed in that paper.

\begin{thm*}
Let $\Sigma$ be a rational pointed fan  in $\RR^d$ and $a \in \ZZ^d$.
Then
\begin{eqnarray*}
H_\mm^i(K[\Sigma])
&\iso&
\bigoplus_{C\in \Sigma}
\bigoplus_{a \in -\relint(C)}
\Tilde H^{i-\dim C-1}(\Delta(\str_\Sigma(C))) \otimes_K K(-a)
\end{eqnarray*}
as $\ZZ^d$-graded $K$-modules.
\end{thm*}

Here $\str_\Sigma(C)=\{ D \in \Sigma : C \subseteq D\}$ is the {\em star}
of $C$, which is also a poset ordered by inclusion.
$\Delta(\str_\Sigma(C))$ is the order complex of the poset $\str_\Sigma(C)\setminus\{C\}$,
which is the simplicial complex where the faces are the chains of $\str_\Sigma(C)\setminus\{C\}$,
and
$\Tilde H^{i-\dim C-1}(\Delta(\str_\Sigma(C)))$ denotes the simplicial cohomology of
the simplicial complex $\str_\Sigma(C)$ with respect to $K$.

In Section \ref{canonicalmodule}
we  compute the $\ZZ^d$-graded dualizing complex for a toric face ring $K[\Sigma]$.
If $K[\Sigma]$ is Cohen--Macaulay
we determine the $\ZZ^d$-graded canonical module of $K[\Sigma]$.
This results are generalizations to those ones known for Stanley--Reisner rings
(see \cite{BH,GR84,GR842}).

Finally, in Section \ref{gorensteinsection}
we characterize Gorenstein toric face rings.
For $C \in \Sigma$
let
$\Tilde\chi_\Sigma (C) =\sum_{i=-1}^{\dim \Delta(\str_\Sigma(C))} (-1)^{i}f_i(\Delta(\str_\Sigma(C)))$
be the reduced Euler-characteristic of $\Delta(\str_\Sigma(C))$.
In Theorem \ref{gorensteinfirst} we show:

\begin{thm*}
Let $\Sigma$ be a rational pointed fan in $\RR^d$. Then  $K[\Sigma]$ is Gorenstein if and only if $K[\Sigma]$ is Cohen--Macaulay and there
exists $\sigma\in\big [\bigcap_{C \in \Sigma \text{ maximal}}C\big ]\cap\ZZ^d$ such that
we have for all $C \in \Sigma$
$$
\Tilde\chi_\Sigma (C)
=
\begin{cases}
(-1)^{\dim \Sigma -\dim C} , &\text{ if } \relint(C)\cap \ZZ^d -\sigma \subset \bigcup_{C \in \Sigma} C \cap \ZZ^d;\\
0, &\text{else}.
\end{cases}
$$
In particular, if $\sigma=0$, then
$\Tilde\chi_\Sigma (C)=(-1)^{\dim \Sigma -\dim C}$ for all $C \in \Sigma$.
\end{thm*}

We call the fan $\Sigma$
an {\em Euler fan} if $\Sigma$ is pure (i.e.\ all maximal cones of $\Sigma$ have the same dimension),
and $\Tilde\chi_\Sigma (C)=(-1)^{\dim \Sigma - \dim C}$ for all $C \in \Sigma$.
Now we can show the following:

\begin{thm*}
Let $\Sigma$ be a rational pointed fan in $\RR^d$.
The following statements are equivalent:
\begin{enumerate}
\item
$K[\Sigma]$ is Cohen--Macaulay and $\Sigma$ is an Euler fan;
\item
$K[\Sigma]$ is
Gorenstein and $\omega_{K[\Sigma]} \cong K[\Sigma]$ as $\ZZ^d$-graded modules.
\end{enumerate}
\end{thm*}

The remaining part of Section \ref{gorensteinsection} is concerned with the question how
to relate a
Gorenstein toric face ring $K[\Sigma]$ such that $\omega_{K[\Sigma]} \cong K[\Sigma](-\sigma)$
for some $0\neq \sigma \in \ZZ^d$
to a Gorenstein toric face ring of an Euler fan.

We are grateful to Prof.\ W. Bruns
for inspiring discussions on the subject of the paper.

%
%
%
\section{Toric face rings}
\label{toricfacering}
In this section we introduce toric face rings and study related ring properties.
In particular, we determine the graded prime spectrum and graded radical ideals.
Let $\Sigma$ be a fan in $\RR^d$.
In the following we consider $\Sigma$ sometimes also as a partially ordered set, ordered by inclusion.
We say that $\Sigma$ is {\em rational} and {\em pointed} respectively,
if all cones in $\Sigma$ are rational and pointed respectively.
Choose a set of affine monoids $\Mcc_\Sigma=\{M_C \subseteq \ZZ^d : C \in \Sigma \}$
such that:
\begin{enumerate}
\item
$\cn(M_C)=C$ for $C \in \Sigma$;
\item
If $C,C' \in \Sigma$, $C' \subseteq C$, then $M_{C'} = M_C \cap C'$.
\end{enumerate}
In this situation we call $\Mcc_{\Sigma}$ an {\em monoidal complex} supported by $\Sigma$.
Note that here we do not require that $M_C$ is normal.
Observe that $\Sigma$ has to be rational in this case.
That $\Sigma$ is pointed is equivalent to the fact that all affine monoids $M_C$ are
{\em positive}, i.e.\ $0$ is the only invertible
element in $M_C$. We set $|\Mcc_{\Sigma}|=\bigcup_{C\in \Sigma} M_C$.
We say that $\Mcc_{\Sigma}$ has a property $E$ if all monoids $M_C$ satisfy $E$.
E.g.\ $\Mcc_{\Sigma}$ is called {\em Cohen--Macaulay} and {\em Gorenstein} respectively,
if all affine monoids $M_C$ are Cohen--Macaulay and Gorenstein respectively.

Let $K$ be a field.
The {\em toric face ring $K[\Mcc_{\Sigma}]$} of $\Mcc_{\Sigma}$ (over $K$) is defined as follows:
As a $K$-vectors space $K[\Mcc_{\Sigma}]$ has one basis element $x^a$ for
every $a \in |\Mcc_{\Sigma}|$.
Multiplication in $K[\Mcc_{\Sigma}]$ is defined by
the following rule:
\begin{displaymath}
  x^a \cdot x^b =
  \begin{cases}
    x^{a+b} & \text{if $a$ and $b$ are elements of a common monoid of $\Mcc_{\Sigma}$,}\\
    0 & \text{otherwise.}
  \end{cases}
\end{displaymath}
$K[\Mcc_{\Sigma}]$ is a $\ZZ^d$-graded ring, and throughout this paper
the attributes 'graded' and 'homogeneous' refer to the  $\ZZ^d$-graduation of $K[\Mcc_{\Sigma}]$
if not  stated otherwise.

We are in particular interested in the following situation.
Choose $\Ncc_{\Sigma}=\{ C \cap \ZZ^d : C \in \Sigma\}$. Then by Gordan's lemma the elements of $\Ncc_{\Sigma}$
are all normal affine monoids and
$K[\Sigma]=K[\Ncc_{\Sigma}]$ is called the {\em toric face ring} of $\Sigma$,
which was defined and studied by Stanley in \cite{ST87}.
\smallskip

Now let $\Mcc_{\Sigma}$ be a monoidal complex supported by a rational pointed fan $\Sigma$ in $\RR^d$.
Naturally $\ZZ^d$-graded prime and radical ideals are induced from $\Sigma$.
On the one hand for $C \in \Sigma$, let
$\pp_C=(x^a \in K[\Mcc_\Sigma] : a \not\in M_C)$.
(Note that $\{0\}$ is also a face of $\Sigma$.)
Then $\pp_C$ is $\ZZ^d$-graded and prime, because
$K[\Mcc_\Sigma]/\pp_C$ is isomorphic to the affine monoid ring $K[M_C]$ associated to the affine monoid $M_C$, which is an integral domain.
On the other hand let $\Sigma' \subseteq \Sigma$ be a subfan,
$\Mcc_{\Sigma'}$ the induced monoidal complex and set
$\qq_{\Sigma'}=(x^a \in K[\Mcc_\Sigma] : a \not\in |\Mcc_{\Sigma'}|)$.
Then
$\qq_{\Sigma'}$ is a $\ZZ^d$-graded radical ideal of $K[\Mcc_\Sigma]$ since
$K[\Mcc_\Sigma]/\qq_{\Sigma'}$ is isomorphic to $K[\Mcc_{\Sigma'}]$, which is reduced.
The next lemma shows that every graded prime and radical ideal is of the presented forms.

\begin{lemma}
\label{primeradical_ideals}
Let $\Mcc_\Sigma$ be a monoidal complex supported by a rational pointed fan $\Sigma$ in $\RR^d$.
\begin{enumerate}
\item
The assignment $C\mapsto \pp_C$
is a bijection between
the set of non-empty cones in $\Sigma$ and the set of $\ZZ^d$-graded prime ideals of $K[\Mcc_\Sigma]$.
\item
The assignment $\Sigma'\mapsto \qq_{\Sigma'}$
is a bijection between
the set of non-empty subfans of $\Sigma$
and the set of $\ZZ^d$-graded radical ideals of $K[\Mcc_\Sigma]$.
\end{enumerate}
In particular, $\mm=\pp_0$ is the unique $\ZZ^d$-graded maximal ideal of $K[\Mcc_\Sigma]$,
which is also maximal in the usual sense.
\end{lemma}
\begin{proof}
The injectivity of the assignments in (i) and (ii) is clear and it remains to prove the surjectivity.

Observe that $0=\bigcap_{C \in \Sigma \text{ maximal}} \pp_C \subset K[\Mcc_\Sigma]$
is an irrendundant primary decomposition of $0$, and thus we have that
the $\pp_C$ for the maximal cones $C$ are exactly the minimal prime ideals of
$K[\Mcc_\Sigma]$. Now let $\pp$ be an arbitrary $\ZZ^d$-graded  prime ideal of $K[\Mcc_\Sigma]$.
Then there is  a minimal $\ZZ^d$-graded  prime ideal $\pp_C$  of $K[\Mcc_\Sigma]$
such that $\pp_C \subseteq \pp$.
Since $\bar \pp$ is a $\ZZ^d$-graded  prime ideal of $K[\Mcc_\Sigma]/\pp_C\cong K[M_C]$,
the conclusion of (i) follows now from \cite[Theorem 6.1.7]{BH}.
The surjectivity of the assignment in (ii) follows from (i), since each radical
ideal is the intersection of its minimal prime ideals.
\end{proof}

Recall that $\Sigma$ is called {\em pure}
if all maximal cones in $\Sigma$  have the same dimension.
Using the irrendundant primary decomposition of $0$
we immediately obtain:
\begin{corollary}
\label{CMpure}
Let $\Mcc_\Sigma$ be a monoidal complex supported by a rational pointed fan $\Sigma$ in $\RR^d$
such that $K[\Mcc_\Sigma]$ is a Cohen--Macaulay ring.
Then $\Sigma$ is pure.
\end{corollary}

In subsequent sections we need the following easy result.

\begin{lemma}
\label{primehelper}
Let $\Mcc_\Sigma$ be a monoidal complex supported by a rational pointed fan $\Sigma$ in $\RR^d$
and $C_1,\dots,C_j,D \in \Sigma$. Then we have in $K[\Mcc_\Sigma]$ the following
equalities:
\begin{enumerate}
\item
$\pp_{C_1\cap\dots\cap C_j}= \pp_{C_1} +\dots+ \pp_{C_j}$.
\item
$\pp_D+\bigcap_{i=1}^t \pp_{C_i}= \bigcap_{i=1}^t (\pp_D +\pp_{C_i})$.
\end{enumerate}
\end{lemma}
\begin{proof}
Observe that the $\pp_{C_i}$ are $\ZZ^d$-graded prime ideals of $K[\Mcc_\Sigma]$.
Thus they are monomial ideals in this ring.
Their bases as $K$-vector spaces are subsets of the set of
monomials $x^a$, $a\in |\Mcc_\Sigma|$. Using this fact it is easy to check the
claimed equalities.
\end{proof}
%
%
%
\section{Shellable fans}
\label{shellfan}
In this section we present combinatorial conditions
which imply that a toric face ring is Cohen--Macaulay and sequentially Cohen--Macaulay respectively.
The results generalize well-known results for Stanley--Reisner rings.

Let $\Sigma$ be a rational pointed fan in $\RR^d$. In the following $\partial C$ denotes the boundary of a cone $C \in \Sigma$.
The fans $\fan(C)$ and $\fan(\partial C)$ are the set of faces of $C$ and $\partial C$ respectively.
Observe that they are subfans of $\Sigma$.
Assume that $\Sigma$ is {\em pure} $k$-dimensional, i.e.\
the facets of $\Sigma$ have all the same dimension $k$.
Recall from \cite{BJWA97}
that a {\em shelling} of $\Sigma$
is a linear ordering $C_1,\dots, C_s$
of the facets of $\Sigma$ such that either $k=0$,
or the following two conditions are satisfied:
\begin{enumerate}
\item
$\fan(\partial C_1)$ has a shelling.
\item
For $1<j\leq s$ the fan $\fan(\partial C_j)$ is pure $k-1$-dimensional and
there exists a shelling $D_{1},\dots, D_{t_j}$ of $\fan(\partial C_j)$
such that
$$
\emptyset
\neq
\bigcup_{i=1}^{j-1} \big [\fan(C_i) \cap \fan(C_j)\big ]
=
\bigcup_{l=1}^{r_j} \fan(D_l)
\text{ for some } 1\leq r_j \leq t_j.
$$

\end{enumerate}
$\Sigma$ is called {\em shellable} if it is pure and has a shelling.
Note that the assumption that $\fan(\partial C_1)$ has a shelling could be omitted,
since by Brugesser--Mani the boundary of a cross-section of $C_1$ has a shelling and
thus $C_1$ has also one. Also $\fan(\partial C)$ is always pure of dimension $\dim C-1$.
So we may only ask that  for $1<j\leq s$
there exists a shelling $D_{1},\dots, D_{t_j}$ of $\fan(\partial C_j)$
such that
$$
\emptyset
\neq
\bigcup_{i=1}^{j-1} \big [\fan(C_i) \cap \fan(C_j)\big ]
=
\bigcup_{l=1}^{r_j} \fan(D_l)
\text{ for some } 1\leq r_j \leq t_j.
$$
\begin{example}
\
\begin{enumerate}
\item
Every rational pointed fan $\Sigma$ in $\RR^d$ of dimension $\leq 1$ is shellable for
trivial reasons. A non shellable fan is given by considering e.g.\ two rational pointed
cones of dimension $2$ which meet only in $\{0\}$ and all their faces.
\item
We already noted that for a single rational pointed cone $C \subseteq \RR^d$,
the fan $\fan(\partial C)$ is shellable. The other extreme case is that
if $\Sigma$ is a {\em simplicial} fan, i.e.\ each cone $C \in \Sigma$
is simplicial. Then $\Sigma$ can be considered as a simplicial complex.
This case reduces then to the notion of shellable simplicial complexes.
See Björner and Wachs \cite{BJWA96, BJWA97} for details and further examples.
\end{enumerate}
\end{example}
As the shellability of a simplicial complex implies that the corresponding
Stanley--Reisner ring is Cohen--Macaulay, we have that the shellability of fans imply
Cohen--Macaulayness of the corresponding toric face rings if some mild assumptions are satisfied.
More precisely, we have:

\begin{theorem}
\label{shellcm}
Let $\Mcc_\Sigma$ be a Cohen-Macaulay monoidal complex supported by a
rational pointed shellable fan $\Sigma$ in $\RR^d$.
Then $K[\Mcc_\Sigma]$ is Cohen--Macaulay.
In particular, the toric face ring $K[\Sigma]$ is Cohen--Macaulay (independent of $\cha K$).
\end{theorem}
\begin{proof}
Let $\Sigma$ be pure $k$-dimensional and $C_1,\dots, C_s$ be a shelling of $\Sigma$.
We may assume that $k>1$.
For $j=1,\dots,s$ consider the subfans $\Sigma_j=\bigcup_{i=1}^j \fan(C_i)$ of $\Sigma$
and the
corresponding submonoidal complexes $\Mcc_j=\Mcc_{\Sigma_j}$ of $\Mcc_\Sigma$.

We show by induction on $j$ that $K[\Mcc_j]$ is Cohen--Macaulay of dimension $k$.
Then the case $j=s$ implies that $K[\Mcc_{\Sigma}]$ is Cohen--Macaulay.
The case $j=1$ is easy, since
$K[\Mcc_1]=K[M_{C_1}]$ is a Cohen--Macaulay ring by assumption.
Furthermore we have $\dim K[M_{C_1}]=\rank M_{C_1}= \dim C_1=k$.

Let $j>1$.
We denote by $\Pi_j$ the fans
$\bigcup_{i=1}^{j-1} \big [\fan(C_i) \cap \fan(C_j)\big ]$
for $j=1,\dots,s$ with
corresponding monoidal Cohen--Macaulay complexes
$\Wcc_j=\{M_C : C \in \Pi_j\}$.
It is important to observe that $\Pi_j$ is pure $k-1$-dimensional, and again shellable
by the very definition of shellability of $\Sigma$.
Observe that
$$
K[\Mcc_j]=K[\Mcc_\Sigma]/ \bigcap_{i=1}^j  \pp_{C_i},\quad
K[M_{C_j}]=K[\Mcc_\Sigma]/ \pp_{C_j}
$$
and
$$
K[\Wcc_j]
=
K[\Mcc_\Sigma]/ \bigcap_{i=1}^{j-1}  (\pp_{C_i} + \pp_{C_j})
=
K[\Mcc_\Sigma]/ (\bigcap_{i=1}^{j-1}  \pp_{C_i}) + \pp_{C_j}
$$
where the latter equality follows from Lemma \ref{primehelper}.
Consider the following homomorphisms of  $K[\Mcc_{\Sigma}]$-modules
$$
\alpha \colon
K[\Mcc_j]
\to
K[\Mcc_{j-1}]
\bigoplus
K[M_{C_j}],\
\overline{a} \mapsto
(\overline{a}, -\overline{a}),
$$
$$
\beta\colon
K[\Mcc_{j-1}] \bigoplus K[M_{C_j}]
\to
K[\Wcc_j],\
(\overline{a}, \overline{b})
\mapsto
(\overline{a + b}),
$$
and the associated short exact sequence:
$$
0
\to
K[\Mcc_j]
\overset{\alpha}{\to}
K[\Mcc_{j-1}] \bigoplus K[M_{C_j}]
\overset{\beta}{\to}
K[\Wcc_j]
\to
0.
$$
By induction we have that
$K[\Mcc_{j-1}]$ is Cohen--Macaulay of dimension $k$.
By assumption $K[M_{C_j}]$
is Cohen--Macaulay of dimension $k$ since $\rank M_{C_j}=k$.
The crucial point is that  $\Pi_j$ is shellable and pure $k-1$-dimension.
Thus we may also conclude by induction on $k$ that $K[\Wcc_j]$ is Cohen--Macaulay of
dimension $k-1$.
Hence standard arguments yield that $K[\Mcc_j]$ is Cohen--Macaulay of dimension $k$.

It remains to note that by Hochster the affine monoid rings
$K[C\cap \ZZ^d]$ are Cohen--Macaulay of dimension $\dim C$ (independent of $\cha K$),
since $C\cap \ZZ^d$ is normal,
to conclude that $K[\Sigma]$ is Cohen--Macaulay if $\Sigma$ is shellable (independent of $\cha K$).
\end{proof}

We saw that the shellability of $\Sigma$ implies that $K[\Sigma]$ is Cohen--Macaulay (over any field).
One can weaken the notion of shellability. As defined by Björner and Wachs \cite{BJWA96},
one calls $\Sigma$ {\em non-pure shellable} if $\Sigma$ admits a shelling
without the assumption that $\Sigma$ has to be pure, i.e.\
$\Sigma$ is non-pure shellable if there
is a linear ordering $C_1,\dots, C_s$
of the facets of $\Sigma$ such that either $\dim \Sigma=0$,
or the following two conditions are satisfied:
\begin{enumerate}
\item
$\fan(\partial C_1)$ has a shelling.
\item
For $1<j\leq s$
there exists a shelling
$D_{1},\dots, D_{t_j}$ of the (pure) fan $\fan(\partial C_j)$
such that
$
\emptyset
\neq
\bigcup_{i=1}^{j-1} \fan(C_i) \cap \fan(C_j)
=
\bigcup_{l=1}^{r_j} \fan(D_l)$
for some $1\leq r_j \leq t_j$.
\end{enumerate}

We need the following Lemma about non-pure shellings which is analogue to
the corresponding results \cite[Lemma 2.6]{BJWA96}
for simplicial complexes.

\begin{lemma}
\label{shellinghelper}
Let $\Sigma$ be a non-pure shellable fan in $\RR^d$
with shelling $C_1,\dots, C_s$.
Then there exists a permutation $\sigma$ of $\{1,\dots,s\}$
such that $C_{\sigma(1)},\dots, C_{\sigma(s)}$ is a shelling of $\Sigma$
and $\dim C_{\sigma(1)} \geq \dots \geq \dim C_{\sigma(s)}$.
\end{lemma}
In other words, if there exists a shelling of $\Sigma$, then
there exists also one such that the dimensions are decreasing.
\begin{proof}
Let $\dim \Sigma =k$.
Let $C_{m_1},\ldots,C_{m_s}$ be
the rearrangement obtained
by taking first all facets of dimension $k$ in the
order of the given shelling, then all facets of dimension $k-1$ and continuing in this way.
Let $\sigma$ be the permutation with $\sigma(j)=m_j$ for $j=1,\dots,s$.

We claim that $C_{m_1},\ldots,C_{m_s}$ is again a shelling of $\Sigma$.
For this it is enough to prove that for all $1< j\leq s$
we have
\begin{equation}
\label{eqn1}
\bigcup_{i=1}^{j-1} \big [\fan(C_{m_i}) \cap \fan(C_{m_j})\big ]
=
\bigcup_{i=1}^{m_j-1} \big [\fan(C_i) \cap \fan(C_{m_j})\big ]
\end{equation}
since the right hand side is the pure shellable subfan
$\bigcup_{l=1}^{r_{m_j}} \fan(D_l)$
of $\fan(\partial C_{m_j})$ by assumption. It is enough to show that the facets
of the left and right fan respectively
are contained in the right and left fan respectively.

The facets of the right hand side of (\ref{eqn1}) are of the form
$C_{i} \cap C_{m_j}=D_l$ for some $i<m_j$ and $1 \leq l \leq r_{m_j}$.
But then $\dim C_i\geq \dim D_l +1 = \dim C_{m_j}$ and thus by the choice of the
of the new order of the facets of $\Sigma$
there exists an $h<j$ such that
$m_h=i$ and $D_l \in \bigcup_{i=1}^{j-1} \big [\fan(C_{m_i}) \cap \fan(C_{m_j})\big ]$.

The facets of the left hand side of (\ref{eqn1}) are of the form
$D=C_{m_i} \cap C_{m_j}$ for some $i<j$. Note that then $\dim C_{m_i} \geq \dim C_{m_j}$.
If $m_i<m_j$, then $D$ is trivially an element of the right hand side
of (\ref{eqn1}).
Assume that there exist $j$ and $i$ with $m_i> m_j$ and $C_{m_j}$ comes after $C_{m_i}$ in the ``new'' shelling.
Then $\dim C_{m_i}> \dim C_{m_j}$.
If
$D\in \bigcup_{l=1}^{r_{m_j}} \fan(D_l)$, then we are done.
Assume now that
$D \not\in \bigcup_{l=1}^{r_{m_j}} \fan(D_l)$.
Hence there exists
$$
j<i \text{ such that } \dim C_{i} > \dim C_{j}
\text{ and }
D=C_{i} \cap C_{j} \not\in
\bigcup_{h=1}^{j-1} \big [\fan(C_h) \cap \fan(C_{j})\big ].
$$
Choose   a pair $j,i$ with $i$ minimal over all such choices.
Since $\Sigma$ is shellable, there exists an $n<i$ such that
$D \subseteq C_n\cap C_i$ and
$\dim C_n\cap C_i=\dim C_i -1$.
Then
$$
C_n\cap C_j
\not\in
\bigcup_{h=1}^{j-1} \big [\fan(C_h) \cap \fan(C_{j})\big ],\
\dim C_n > \dim C_j
\text{ and thus also } n >j.
$$
This is a contradiction to the minimality of $i$.
Thus we conclude that we found a new shelling of $\Sigma$
such that the dimension of the facets in the new order are decreasing.
\end{proof}

On the algebra side
one defines  the notion of
sequentially  Cohen--Macaulay modules.
Let $K$ be a field and $S$ be a finitely generated $\ZZ^d$-graded $K$-algebra.
Let $M$ be a finitely generated  $\ZZ^d$-graded $S$-module.
A finite filtration
$$
0=M_0 \subset M_1 \subset \dots \subset M_r=M
$$
of $\ZZ^d$-graded submodules of $M$ is called a {\em CM-filtration},
if each quotient $M_i/M_{i-1}$ is Cohen--Macaulay and
$\dim(M_1/M_0)<\dim( M_2/M_1) <\dots<\dim (M_r/M_{r-1})$.
The module $M$ is called {\em sequentially Cohen--Macaulay},
if $M$ has a CM-filtration.
In the definition of a sequentially Cohen--Macaulay module
we required to have strict inequalities of the dimensions.
A little weaker assumption implies also sequentially Cohen--Macaulayness,
as the next result shows.

\begin{lemma}
\label{seqcmhelper}
Let $N_0 \subset N_1 \subset \dots \subset N_s$ be finitely generated $\ZZ^d$-graded
$S$-modules and let $k$ be an integer  such that all $N_i/N_{i-1}$ are Cohen--Macaulay of dimension $k$.
Then $N_s/N_0$ is Cohen--Macaulay of dimension $k$.

In particular, if a finitely generated $\ZZ^d$-graded $S$-module $M$ has a
filtration
$0=M_0 \subset M_1 \subset \dots \subset M_r=M$
of $\ZZ^d$-graded submodules of $M$ such that each quotient $M_i/M_{i-1}$ is Cohen--Macaulay
and
$\dim(M_1/M_0)\leq \dim( M_2/M_1) \leq \dots \leq\dim (M_r/M_{r-1})$, then $M$ is
sequentially Cohen--Macaulay.
\end{lemma}
\begin{proof}
We prove by induction on $i \in \{1,\dots, s\}$ that
$N_i/N_{0}$ is Cohen--Macaulay of dimension $k$.

For $i=1$ we have that $N_1/N_0$ is Cohen--Macaulay of dimension $k$
by assumption. Let $i>1$ and
consider the following short exact sequence
$$
0
\to
N_{i-1}/ N_0
\to
N_i/ N_0
\to
N_{i}/ N_{i-1}
\to
0.
$$
By the induction hypothesis
$N_{i-1}/ N_0$ is Cohen--Macaulay of dimension $k$.
By assumption
$N_{i}/ N_{i-1}$ is Cohen--Macaulay of dimension $k$.
Thus standard arguments yield that
$N_{i}/ N_0$ is Cohen--Macaulay of dimension $k$.

Assume now that
a finitely generated $\ZZ^d$-graded $S$-module $M$ has a
filtration
$0=M_0 \subset M_1 \subset \dots \subset M_r=M$
of $\ZZ^d$-graded submodules of $M$ such that each quotient $M_i/M_{i-1}$ is Cohen--Macaulay
and
$\dim(M_1/M_0)\leq \dim( M_2/M_1) \leq \dots \leq\dim (M_r/M_{r-1})$.
Let $i_1<\dots< i_s$ be those numbers such that
$\dim(M_{i_j}/M_{i_j - 1})< \dim(M_{i_j+1}/M_{i_j})$.
Then we know by what we proved above
that $M_{i_{j}}/M_{i_{j-1}}$ is Cohen--Macaulay of dimension
$\dim(M_{i_j}/M_{i_j - 1})$
and
$$
0= M_0 \subset M_{i_1} \subset \dots \subset M_{i_s}=M
$$
is a CM-filtration of $M$, which shows that $M$ is sequentially Cohen--Macaulay.
\end{proof}

Now we apply our results to toric face rings.

\begin{theorem}
\label{shellseqcm}
Let $\Mcc_\Sigma$ be a Cohen--Macaulay monoidal complex supported by a rational pointed
non-pure shellable fan $\Sigma$ in $\RR^d$.
Then $K[\Mcc_\Sigma]$ is sequentially Cohen--Macaulay.
In particular, the toric face ring $K[\Sigma]$ is sequentially  Cohen--Macaulay
(independent of $\cha K$).
\end{theorem}

Observe that for $M_C \in \Mcc_\Sigma$ the ring $K[M_C]$ is sequentially Cohen--Macaulay if and only
if $K[M_C]$ is Cohen--Macaulay, because this ring is an integral domain.
Hence one can not weaken the assumption Cohen--Macaulay on the monoidal complex
in this direction.

\begin{proof}
Let $\Sigma$ be $k$-dimensional and $C_1,\dots, C_s$ be a non-pure shelling of $\Sigma$.
By Lemma \ref{shellinghelper} we may assume  $\dim C_1 \geq \dots \geq \dim C_s$.
The cases $k\leq 1$ are trivial, thus let $k>1$.
For $j=1,\dots,s$ consider again
the subfans $\Sigma_j=\bigcup_{i=1}^j \fan(C_i)$ of $\Sigma$
and the corresponding submonoidal complexes $\Mcc_j=\Mcc_{\Sigma_j}$ of $\Mcc_\Sigma$.

Now we show by induction on $j$ that $K[\Mcc_j]$ is sequentially Cohen--Macaulay of dimension $k$.
Then the case $j=s$ implies that $K[\Mcc_{\Sigma}]$ is sequentially Cohen--Macaulay.
Since $K[\Mcc_1]=K[M_{C_1}]$ is a Cohen--Macaulay ring by assumption  we proved the case $j=1$.

Let $j>1$.
Consider the fans
$\Pi_j=\bigcup_{i=1}^{j-1} \big [\fan(C_i) \cap \fan(C_j)\big ]$
for $j=1,\dots,s$ with
corresponding monoidal Cohen--Macaulay complexes
$\Wcc_j=\{M_C : C \in \Pi_j\}$.
It is important to observe that $\Pi_j$ is pure of dimension $\dim C_j-1$ and shellable.
Thus we know by Theorem  \ref{shellcm} that $K[\Wcc_j]$ is Cohen--Macaulay of dimension $\dim C_j-1$.
We consider the filtration
$$
0=
\bigcap_{i=1}^s \pp_{C_i}
\subset
\bigcap_{i=1}^{s-1} \pp_{C_i}
\subset
\dots
\subset
\bigcap_{i=1}^2 \pp_{C_i}
\subset
\pp_{C_1}
\subset
K[\Mcc_{\Sigma}]
$$
and claim that
$\bigcap_{i=1}^{j-1} \pp_{C_i} / \bigcap_{i=1}^{j} \pp_{C_i} $ is Cohen--Macaulay
of dimension $\dim C_j$. Then it follows from Lemma \ref{seqcmhelper} that
$K[\Mcc_{\Sigma}]$ is sequentially Cohen--Macaulay.
Since
$K[\Mcc_{\Sigma}]/\pp_{C_1}$ is Cohen--Macaulay of dimension $\dim C_1$
by the assumption that $\Mcc_{\Sigma}$ is a Cohen--Macaulay monoidal complex,
we may assume that $j \geq 2$.
Observe that
$$
\bigcap_{i=1}^{j-1} \pp_{C_i} / \bigcap_{i=1}^{j} \pp_{C_i}
\cong
\bigl ( \bigcap_{i=1}^{j-1} \pp_{C_i} \bigr ) + \pp_{C_j} / \pp_{C_j}
\cong
\bigl (\bigcap_{i=1}^{j-1} (\pp_{C_i} + \pp_{C_j}) \bigr ) / \pp_{C_j}
$$
where the last isomorphism follows from
Lemma \ref{primehelper}.
Consider the following short exact sequence
$$
0
\to
\bigl (\bigcap_{i=1}^{j-1} (\pp_{C_i} + \pp_{C_j}) \bigr ) / \pp_{C_j}
\to
K[\Mcc_{\Sigma}] / \pp_{C_j}
\to
K[\Mcc_{\Sigma}]/ \bigl (\bigcap_{i=1}^{j-1} (\pp_{C_i} + \pp_{C_j}) \bigr )
\to 0.
$$
Now
$K[\Mcc_{\Sigma}] / \pp_{C_j}=K[M_{C_j}]$ is Cohen--Macaulay
of dimension $\dim C_j$ by assumption.
Furthermore,
$K[\Mcc_{\Sigma}]/ \bigl (\bigcap_{i=1}^{j-1} (\pp_{C_i} + \pp_{C_j}) \bigr )=K[\Pi_j]$
is Cohen--Macaulay of dimension $\dim C_j-1$ as was noted above.
By standard arguments
$\bigl (\bigcap_{i=1}^{j-1} (\pp_{C_i} + \pp_{C_j}) \bigr ) / \pp_{C_j}$
is Cohen--Macaulay of dimension $\dim C_j$.
\end{proof}

%
%
%
\section{Local cohomology}
\label{localcoho}

In this section we define a complex which computes the local cohomology of a toric face ring
in a similar way as it was done for an affine monoid ring in \cite[Section 6.2]{BH}.
Let  $\Sigma$ be a nontrivial rational pointed fan in $\RR^d$.
We consider the intersection $X$ of $\Sigma$
with the unit sphere $S^{d-1}$ and the set
$\Gamma_\Sigma=\big\{\relint(C)\cap S^{d-1} : C\in\Sigma\big\}.$
Here $\relint(C)$ denotes the relative interior of $C$ with respect to the
subspace topology on the vector space generated by $C$.
For a cone $C\in \Sigma$ we denote by $e_C=\relint(C)\cap S^{d-1}$ the corresponding element of $\Gamma_\Sigma$.
Set
$\Gamma_\Sigma^i=\big\{e\in\Gamma_\Sigma:\bar e \text{ homeomorph with }B^i\big\}$
where $B^i$ is the $i$-dimensional ball in $\RR^i$. The elements $e_C$ are called {\em open cells}.
Then $(X,\Gamma_\Sigma)$ is a {\em finite regular cell complex}.
The dimension of $\Gamma_\Sigma$ is given by
$\dim \Gamma_\Sigma=\max\{i:\Gamma_\Sigma^i\not=\emptyset\}(= \dim \Sigma -1)$.
An element $e_{C'}$ is called a {\em face} of $e_C$ if $e_{C'} \subset \bar e_C$, i.e.\
$C'$ is a face of $C$.

There exists an {\em incidence function} $\varepsilon$ on $\Gamma_\Sigma$,
i.e.\ $\varepsilon$ assigns to each pair $(e_C,e_{C'})$ with
$e_C \in \Gamma_\Sigma^i$ and $e_{C'} \in \Gamma_\Sigma^{i-1}$ for some $i\geq 0$ a number
$\varepsilon(e_C,e_{C'}) \in \{0,\pm 1\}$, such that the following is satisfied:
$\varepsilon(e_C,e_{C'}) \neq 0$ if and only if $e_{C'}$ is a face of $e_C$,
$\varepsilon(e_C,\emptyset)=1$ for all $0$-cells $e_C$,
and if
$e_C \in \Gamma_\Sigma^i$, $e_{C'} \in \Gamma_\Sigma^{i-2}$, then
$$
\varepsilon(e_C,e_{C_1})\varepsilon(e_{C_1},e_{C'})
+
\varepsilon(e_C,e_{C_2})\varepsilon(e_{C_2},e_{C'})
=0
$$
where $e_{C_1}, e_{C_2}$ are those uniquely determined  $(i-1)$-cells such that $e_{C_j}$ is a face of $e_{C}$
and $e_{C'}$ is a face of $e_{C_j}$.
Given a cell complex $(X,\Gamma_\Sigma)$ of dimension $k-1$, we define the augmented oriented chain complex of $\Gamma_\Sigma$ by
$$
\CD
\mathcal{C}_{\lpnt}(\Gamma_\Sigma)
\colon
0@>>>\mathcal{C}_{k-1}@>\partial>>\mathcal{C}_{k-2}@>>>\cdots@>>>\mathcal{C}_{0}@>\partial>>\mathcal{C}_{-1}@>>>0
\endCD
$$
where we set
$
\mathcal{C}_{i}=\bigoplus_{e_C\in \Gamma_\Sigma^i}\ZZ e_C$
for $i=0,\ldots,k-1$, $\mathcal{C}_{-1}=\ZZ$,
and for
$e_{C}\in \Gamma_\Sigma^i$ the differential is given by
$\partial(e_{C})=\sum_{e_{C'}\in\Gamma_\Sigma^{i-1}}\varepsilon(e_{C},e_{C'})e_{C'}$.
We set $\Tilde H_i(\Gamma_\Sigma)=\Tilde H_i(\mathcal{C}_{\lpnt}(\Gamma_\Sigma))$
for the homology of $\mathcal{C}_{\lpnt}(\Gamma_\Sigma)$.
The next corollary will be useful in the following.
\begin{lemma}[Corollary 6.2.4. \cite{BH}]
\label{exact}
If $\Sigma$ is a fan corresponding to one cone, then
$\Tilde H_i(\Gamma_\Sigma)=0$ for every $i\geq -1$, i.e.\
$\mathcal{C}_{\lpnt}(\Gamma_\Sigma)$ is exact.
\end{lemma}

Let $\Mcc_\Sigma$ be a monoidal complex supported by $\Sigma$ and $K$ a field.
For $C\in\Sigma$ we denote by $K[\Mcc_\Sigma]_C$
the homogeneous localization $K[\Mcc_\Sigma]_{(\pp_C)}$. This is
the localization with respect to the set of homogeneous elements of $K[\Mcc_\Sigma]$ not belonging to $\pp_C$.
Let
$$
L^t=\bigoplus_{e_C\in \Gamma_\Sigma^{t-1}}K[\Mcc_\Sigma]_C,\phantom{aaaaaa} t=0,\ldots,k,
$$
and define $\partial:L^{t-1}\to L^t$ by specifying its component $\partial_{C',C}:K[\Mcc_\Sigma]_{C'}\to K[\Mcc_\Sigma]_C$ to be
$\partial_{C',C}=\varepsilon(C,C')\nat,$
where $\varepsilon$ is the incidence function chosen for $\Sigma$.
It follows from the properties of  an incidence function,
that
$$
\CD
L^{\pnt}(\Mcc_\Sigma): 0@>>> L^0 @>\partial>> L^1 @>>>\cdots @>>> L^{k-1} @>\partial>> L^k @>>> 0
\endCD
$$
is  a complex. We set $L^{\pnt}(\Sigma)=L^{\pnt}(\Ncc_\Sigma)$ for the monoidal complex
$\Ncc_{\Sigma}=\{ C \cap \ZZ^d : C \in \Sigma\}$ associated to $\Sigma$.
By local cohomology groups of $K[\Mcc_\Sigma]$
we always mean the local cohomology with respect to the
maximal $\ZZ^d$-graded ideal $\mm$  of $K[\Mcc_\Sigma]$.
Observing that with minor modification the proof of \cite[Theorem 6.2.5]{BH}
works also for toric face rings instead of affine monoid rings, we get:

\begin{theorem}
\label{mainlocalcoho}
Let $\Mcc_\Sigma$ be a monoidal complex supported by a rational pointed fan $\Sigma$ in $\RR^d$.
Then we have for every $K[\Mcc_\Sigma]$-module $M$, and all $i\ge 0$ that
$$
H_\mm^i(M) \cong H^i(L^{\pnt}(\Mcc_\Sigma)\otimes_{K[\Mcc_\Sigma]} M).
$$
\end{theorem}

\begin{proof}
We follow the pattern of the proof of \cite[Theorem 6.2.5]{BH}.
We only sketch the arguments and the main differences to the latter proof.
In order to prove the isomorphism,
we show that the functors $H^i(L^{\pnt}(\Mcc_\Sigma)\otimes_{K[\Mcc_\Sigma]}\_)$
are the right derived functors of $H_\mm^0(\_)$.
For this we have to show firstly that $H_\mm^0(M)=H^0(L^{\bullet}\otimes_{K[\Mcc_\Sigma]} M)$,
secondly that by tensoring a short exact sequence
$0\to M_1 \to M_2\to M_3\to 0$ of $K[\Mcc_\Sigma]$-modules with $L^{\pnt}(\Mcc_\Sigma)$,
we get an induced long exact homology sequence,
and thirdly that $H^i(L^{\pnt}(\Mcc_\Sigma)\otimes M)=0$
for all integers $i>0$ if M is an injective $K[\Mcc_\Sigma]$-module.
The first and the second part may be done exactly as in \cite{BH}.
The only notable difference appears in the third part.

Since every injective module over a Noetherian ring is a direct sum of indecomposable injective modules,
it suffices to consider the indecomposable modules
$E(K[\Mcc_\Sigma]/\pp)$, where $\pp$ is a prime ideal of $K[\Sigma]$.
Let $y \in K[\Mcc_\Sigma]$.
Since $\Ass E(K[\Mcc_\Sigma]/\pp)=\{\pp\}$,
every element of $E(K[\Mcc_\Sigma]/\pp)$ is annihilated by some power of
$y$ if $y\in\pp$. If
$y\not\in\pp$, then
multiplication by $y$
on $E(K[\Mcc_\Sigma]/\pp)$ is an isomorphism.
(It is a injective since $y\not\in\pp$ and surjective since $E(K[\Mcc_\Sigma]/\pp)$ is indecomposable.)
Denote by $\pp^*$ the ideal generated by the set of homogeneous elements belonging to $\pp$. Then
$$
K[\Mcc_\Sigma]_C\otimes_{K[\Mcc_\Sigma]} E(K[\Mcc_\Sigma]/\pp)=\begin{cases}
E(K[\Mcc_\Sigma]/\pp), &\text{ if } \pp^*\subset \pp_C,\\
0, &\text{ if } \pp^*\not\subset \pp_C.
\end{cases}
$$
Since $\pp^*$ is a graded prime ideal,
it follows from Lemma \ref{primeradical_ideals} that
we obtain a cone $C_{\pp^*} \in \Sigma$ such that
$\pp^*=\pp_{C_{\pp^*}}$. Thus
$$
K[\Mcc_\Sigma]_C\otimes E(K[\Mcc_\Sigma]/\pp)=\begin{cases}
E(K[\Mcc_\Sigma]/\pp), &\text{ if } C\subset C_{\pp^*},\\
0, &\text{ if }  C\not\subset C_{\pp^*}.
\end{cases}
$$
Hence we can assume that $\Sigma$ is the fan associated to the cone $C_{\pp^*}$,
and this case was treated in the proof of \cite[Theorem 6.2.5]{BH}.
\end{proof}

Next we are interested in the following question:
Let  $\Sigma$ be a rational pointed fan in $\RR^d$ and
suppose that $\Sigma$ is the union of two subfans, $\Sigma =\Sigma_1 \cup\Sigma_2$.
Then one is interested in the relationship between
the local cohomology of toric face rings with respect to $\Sigma$ and
the local cohomology of toric face rings with respect to $\Sigma_1,\Sigma_2$ and $\Sigma_1\cap\Sigma_2$,
i.e.\ we are searching for a Mayer-Vietoris type of formulae.

Note the following. Let $\Sigma'$ be an arbitrary subfan of $\Sigma$. Then $\Mcc_{\Sigma'}$
denotes the induced submonoidal complex associated to $\Sigma'$.
Observe that the toric face ring $K[\Mcc_{\Sigma'}]$ is a residue class ring of $K[\Mcc_{\Sigma}]$,
and the residue class  of $\mm$ is the maximal ideal of $K[\Mcc_{\Sigma'}]$.
Thus the local cohomology groups of $K[\Mcc_{\Sigma'}]$ coincide with
$H^i_\mm(K[\Mcc_{\Sigma'}])$.
Because of this fact and to avoid cumbersome notation we always
write $H^i_\mm(K[\Mcc_{\Sigma'}])$  for the local
cohomology of $K[\Mcc_{\Sigma'}]$.

\begin{proposition}[The Mayer-Vietoris Sequence]
\label{MV}
Let $\Mcc_\Sigma$ be a monoidal complex supported by a rational pointed fan $\Sigma$ in $\RR^d$.
Suppose that $\Sigma$ is the union of two subfans, $\Sigma =\Sigma_1 \cup\Sigma_2$.
Then there is an exact sequence of $\ZZ^d$-graded $K[\Mcc_\Sigma]$-modules
\begin{align*}
\cdots\to
H_{\mm}^{i-1}(K[\Mcc_{\Sigma_1\cap\Sigma_2}])
\to
H_{\mm}^{i}(K[\Mcc_\Sigma]) &
\to
H_{\mm}^{i}(K[\Mcc_{\Sigma_1}])\oplus H_{\mm}^{i}(K[\Mcc_{\Sigma_2}])\\
&\to
H_{\mm}^{i}(K[\Mcc_{\Sigma_1\cap\Sigma_2}])
\to H_{\mm}^{i+1}(K[\Mcc_\Sigma])\to\cdots
\end{align*}
\end{proposition}

\begin{proof}
Consider for $K[\Mcc_\Sigma]$ the $\ZZ^d$-graded radical ideals
$\qq_{\Sigma_1}$, $\qq_{\Sigma_2}$ and $\qq_{\Sigma_1 \cap \Sigma_2}$
from Lemma \ref{primeradical_ideals}.
Observe that $\qq_{\Sigma_1} \cap \qq_{\Sigma_2}=(0)$ and
$\qq_{\Sigma_1} + \qq_{\Sigma_2}=\qq_{\Sigma_1 \cap \Sigma_2}$. Thus
$$
K[\Mcc_{\Sigma}]=K[\Mcc_\Sigma]/  \qq_{\Sigma_1} \cap \qq_{\Sigma_2} ,\quad
K[\Mcc_{\Sigma_1}]=K[\Mcc_\Sigma]/  \qq_{\Sigma_1},\quad
K[\Mcc_{\Sigma_2}]=K[\Mcc_\Sigma]/  \qq_{\Sigma_2},
$$
and
$$
K[\Mcc_{\Sigma_1 \cap\Sigma_2}]
=
K[\Mcc_\Sigma]/ \qq_{\Sigma_1 \cap \Sigma_2}
=
K[\Mcc_\Sigma]/ \qq_{\Sigma_1} + \qq_{\Sigma_2}.
$$
Consider the following homomorphisms of  $K[\Mcc_{\Sigma}]$-modules
$$
\alpha \colon
K[\Mcc_{\Sigma}]
\to
K[\Mcc_{\Sigma_1 }]
\oplus
K[\Mcc_{\Sigma_2}],\
\overline{a} \mapsto
(\overline{a}, -\overline{a}),
$$
$$
\beta\colon
K[\Mcc_{\Sigma_1 }]
\oplus
K[\Mcc_{\Sigma_2}]
\to
K[\Mcc_{\Sigma_1 \cap \Sigma_2}],\
(\overline{a}, \overline{b})
\mapsto
(\overline{a + b}).
$$
We get the induced short exact sequence of $\ZZ^d$-graded $K[\Mcc_\Sigma]$-modules
$$
0\to
K[\Mcc_{\Sigma}]
\to
K[\Mcc_{\Sigma_1 }]
\oplus
K[\Mcc_{\Sigma_2}]
\to
K[\Mcc_{\Sigma_1 \cap \Sigma_2}]
\to 0
$$
which induce the desired long exact sequence of $\ZZ^d$-graded $K[\Mcc_\Sigma]$-modules
\begin{align*}
\cdots
\to
H_{\mm}^{i-1}(K[\Mcc_{\Sigma_1\cap\Sigma_2}])
\to
H_{\mm}^{i}(K[\Mcc_\Sigma])
&\to
H_{\mm}^{i}(K[\Mcc_{\Sigma_1}])\oplus H_{\mm}^{i}(K[\Mcc_{\Sigma_2}])\\
&\to
H_{\mm}^{i}(K[\Mcc_{\Sigma_1\cap\Sigma_2}])
\to
H_{\mm}^{i+1}(K[\Mcc_\Sigma])
\to\cdots.
\end{align*}
\end{proof}

Recall that we called the toric face ring of the monoidal complex
$\Ncc_{\Sigma}=\{ C \cap \ZZ^d : C \in \Sigma\}$ the toric face ring $K[\Sigma]$ of $\Sigma$.
We also introduced the notation $|\Ncc_{\Sigma}|=\bigcup_{C\in \Sigma} C \cap \ZZ^d$
for the integral points in $\Sigma$.
In the remaining part of this section we concentrate on the local cohomology of $K[\Sigma]$.
We denote by $-\Sigma$ the fan $\{-C:C\in\Sigma\}$.
Note that $H_\mm^i(K[\Sigma])$ is naturally $\ZZ^d$-graded.
At first we present a vanishing result of the local cohomology groups.

\begin{proposition}
\label{vanishinghelper}
Let $\Sigma$ be a rational pointed fan in $\RR^d$ and $a \in \ZZ^d$.
If $a \not\in  |\Ncc_{-\Sigma}|$, then $H_\mm^i(K[\Sigma])_a=0$.
\end{proposition}
\begin{proof}
The case $\dim \Sigma=0$ is trivial, since then $K[\Sigma]=K$.
Assume now that $\dim \Sigma>0$.
Let $C\in \Sigma$ be a cone of maximal dimension $k=\dim \Sigma$ and $a\not\in |\Ncc_{-\Sigma}|$.
Let $\Sigma_1=\Sigma\setminus \{C\}$ and $\Sigma_2=\fan(C)$ be the face poset of $C$.
By Proposition \ref{MV} we have the exact sequence
$$
\cdots \to
H_{\mm}^{i-1}(K[\Sigma_1\cap\Sigma_2])_a
\to
H_{\mm}^{i}(K[\Sigma])_a
\to
H_{\mm}^{i}(K[\Sigma_1])_a
\oplus H_{\mm}^{i}(K[\Sigma_2])_a
\to \cdots
$$
Since $a\not\in - C$ it follows from \cite[Theorem 6.3.4]{BH} that we have
$H_\mm^i(K[\Sigma_2])_a=0$.
Moreover, $a\not\in |\Ncc_{-\Sigma_1}|$ and $a\not\in |\Ncc_{-\Sigma_1\cap \Sigma_2}|$.
Thus by induction on the dimension and another induction on the number of maximal cones in a fan we may assume that
$H_{\mm}^{i-1}(K[\Sigma_1\cap\Sigma_2])_a=0$ and
$H_{\mm}^{i}(K[\Sigma_1])_a=0$.
Hence $H_{\mm}^{i}(K[\Sigma])_a=0$ as desired.
\end{proof}

For $C \in \Sigma$ let $\str_\Sigma(C)=\{D \in \Sigma : C \subseteq D \}$
be the {\em star of $C$}.
Moreover, we set $\Sigma(C)=\Sigma\setminus \str_\Sigma(C)$.
Observe that this a subfan of $\Sigma$.
Note that $\mathcal{C}_{\lpnt}(\Gamma_{ \Sigma(C)})$ is a subcomplex of $\mathcal{C}_{\lpnt}(\Gamma_{\Sigma})$.
We denote by $\mathcal{C}_{\lpnt}(\Gamma_{\str_\Sigma(C)})$  the
factor complex  $\mathcal{C}_{\lpnt}(\Gamma_{\Sigma})/\mathcal{C}_{\lpnt}(\Gamma_{ \Sigma(C)})$
and by
$\mathcal{C}^{\lpnt}(\Gamma_{\str_\Sigma(C)})$ the complex $\Hom_\ZZ(\mathcal{C}_{\lpnt}(\Gamma_{\str_\Sigma(C)}),K)$
. The corresponding homology and cohomology respectively
is denoted by
$\Tilde H_i(\Gamma_{\str_\Sigma(C)})$ and
$\Tilde H^i(\Gamma_{\str_\Sigma(C)})$ respectively.
Observe that for $a\in \relint(C)\cap \ZZ^d$ we have
$\str_\Sigma(C)=\{D \in \Sigma : a \in D\}$.
Hence it makes sense to define slightly more generally for $a\in \ZZ^d$ the {\em star of $a$} as
$\str_\Sigma(a)=\{D\in\Sigma:a\in D\}$. Analogously we define $\Sigma(a)$ and the complexes
$\mathcal{C}_{\lpnt}(\Gamma_{\Sigma(a)})$, $\mathcal{C}_{\lpnt}(\Gamma_{\str_\Sigma(a)})$.
Note that $a\not\in  |\Ncc_{\Sigma(a)}|$. By $\mathcal{C}_{\lpnt}(\Gamma_{\str_\Sigma(-a)})[-1]$
we denote the complex $\mathcal{C}_{\lpnt}(\Gamma_{\str_\Sigma(-a)})$ left shifted with one position
(that is $\mathcal{C}_{i}(\Gamma_{\str_\Sigma(-a)})[-1]=\mathcal{C}_{i-1}(\Gamma_{\str_\Sigma(-a)})$).

\begin{theorem}
\label{thmlcx}
Let $\Sigma$ be a rational pointed fan  in $\RR^d$ and $a \in \ZZ^d$.
Then
$$
L^{\pnt}(\Sigma)_a
\iso
L^{\pnt}( \Sigma(-a) )_a
\oplus
\Hom_\ZZ(\mathcal{C}_{\lpnt}(\Gamma_{\str_\Sigma(-a)})[-1],K)\otimes_K K(-a)
$$
as complexes of $\ZZ^d$-graded $K$-vector spaces.
In particular,
$$
H_\mm^i(K[\Sigma])_a
\iso \Tilde H^{i-1}(\Gamma_{\str_\Sigma(-a)}) \otimes_K K(-a).
$$
\end{theorem}

\begin{proof}
If $\str_\Sigma(-a)=\emptyset$ or $a=0$, then the assertion is trivial.
Otherwise we may assume $\Sigma(-a)$ is a nontrivial subfan of $\Sigma$
and
$\qq_{\Sigma(-a)}$ is a proper ideal of $K[\Sigma]$
with $x^{-a} \in \qq_{\Sigma(-a)} \subset K[\Sigma]$.
At first we show the isomorphism on the module level
for every
$$
L^{t}(\Sigma)_a
=
\bigoplus_{C \in \Sigma,\ \dim C=t} (K[\Sigma]_C)_a
=
\bigoplus_{C \in \Sigma(-a),\ \dim C=t} (K[\Sigma]_C)_a
\oplus
\bigoplus_{C \in \str_\Sigma(-a),\ \dim C=t} (K[\Sigma]_C)_a.
$$
Observe that $\dim_K(K[\Sigma]_C)_a\le 1$.
Indeed, given two elements of $(K[\Sigma]_C)_a$, multiplication by a common
denominator of degree $b$ yields two homogeneous elements of $(K[\Sigma])_{a+b}$. Since $\dim_K(K[\Sigma])_{a+b}\le 1$,
they must be linearly dependent over $K$.

The exact sequence of graded $K[\Sigma]$-modules
$$
\CD
0@>>>\qq_{\Sigma(-a)}@>>>K[\Sigma]@>>>K[\Sigma(-a)]@>>>0
\endCD
$$
induces the exact sequence
$$
\CD
0@>>>({\qq_{\Sigma(-a)}}_C)_a@>>>(K[\Sigma]_C)_a@>>>(K[\Sigma(-a)]_C)_a@>>>0
\endCD
$$
for each $C\in \Sigma$
where ${\qq_{\Sigma(-a)}}_C$ is the localization of $\qq_{\Sigma(-a)}$ in $C$.
We show that $({\qq_{\Sigma(-a)}}_C)_a\not= 0$ if and only if $C\in \str_\Sigma(-a)$.
Since $\dim_K(K[\Sigma]_C)_a\le 1$, we conclude that then
$C\in \str_\Sigma(-a)$ implies $(K[\Sigma(-a)]_C)_a=0$, and that  for $C\in \Sigma(-a)$
we have that $(K[\Sigma]_C)_a \cong (K[ \Sigma(-a)]_C)_a$.
Since
$$\bigoplus_{C \in \str_\Sigma(-a),\ \dim C=t} (K[\Sigma]_C)_a\iso
\Hom_\ZZ(\mathcal{C}_{t-1}(\Gamma_{\str_\Sigma(-a)}),K)
\otimes_K K(-a)
$$
as $K$-vector spaces, the isomorphism on the module level follows.

If $C\in \str_\Sigma(-a)$, then
$({\qq_{\Sigma(-a)}}_C)_a\not=0$,
because $x^{-a}\in \qq_{\Sigma(-a)}$ becomes invertible in ${\qq_{\Sigma(-a)}}_C$
and the inverse has degree $a$.

If $C\in \Sigma(-a)$ we have to show that $({\qq_{\Sigma(-a)}}_C)_a=0$.
Since the multiplicative system  of homogeneous elements not in $(\pp_C)$
does not contain elements of $\qq_{\Sigma(-a)}$,
the ideal ${\qq_{\Sigma(-a)}}_C$ contains no units. Thus ${\qq_{\Sigma(-a)}}_C\not= K[\Sigma]_C$.
In particular, $({\qq_{\Sigma(-a)}}_C)_0=0$ as one easily verifies.
Since multiplication by $x^{-a}$ is injective on ${\qq_{\Sigma(-a)}}$,
multiplication by $\frac{x^{-a}}{1}$ is injective on ${\qq_{\Sigma(-a)}}_C$. On the other hand
$\frac{x^{-a}}{1}({\qq_{\Sigma(-a)}}_C)_a\subset ({\qq_{\Sigma(-a)}}_C)_0=0$.
We conclude that $({\qq_{\Sigma(-a)}}_C)_a=0$.

The final step is to show the isomorphism on the level of homomorphisms,
and hence we really get isomorphism of complexes of $\ZZ^d$-graded $K$-vector spaces.
Consider  $C'\subset C$, $C'\in \Gamma^i_\Sigma$ and $C\in \Gamma^{i+1}_\Sigma$. There are three possibilities:
$C', C\in \Sigma(-a)$, $C', C\in \str_\Sigma(-a)$ or $C'\in \Sigma(-a), C\in \str_\Sigma(-a)$.
In the first case it is easy to see
that the restriction of the differential
$\partial_a:(K[\Sigma]_{C'})_a\to (K[\Sigma]_C)_a$ in $L^{\pnt}(\Sigma)$ is the same as the restriction
of the differential $\partial_a:(K[\Sigma(-a)]_{C'})_a\to (K[\Sigma(-a)]_C)_a$ in $L^{\pnt}(\Sigma(-a))$.
In the second case we identify $\partial_a$ with the corresponding restriction of the differential of $\Hom_\ZZ(\mathcal{C}_{t-1}(\Gamma_{\str_\Sigma(-a)}),K)$.
In the third case we have to show that $\partial_a=0$. We claim that  $(\nat)_a:(K[\Sigma]_{C'})_a\to (K[\Sigma]_C)_a$ is $0$.
Since $C'\subset C\in \str_\Sigma(-a)$, there is some maximal cone which contains both $C'$ and $-a$.
For $b\in C'\cap \ZZ^d$ we deduce $x^{-a}x^b\not=0$. Then
$\frac{x^{-a}}{1}\not=0$ in $K[\Sigma]_{C'}$. On the other hand, the ideal generated
by $x^{-a}$ is contained in $\qq_{\Sigma(-a)}$, and as we have seen above ${\qq_{\Sigma(-a)}}_{C'}\not= K[\Sigma]_{C'}$.
We deduce that $\frac{x^{-a}}{1}$ is not invertible in $K[\Sigma]_{C'}$.
If $y\in (K[\Sigma]_{C'})_a$, then $y\frac{x^{-a}}{1}\in (K[\Sigma]_{C'})_0$. Since $y\frac{x^{-a}}{1}$
is not invertible, it follows that $y\frac{x^{-a}}{1}=0$. Then $\nat(y)\nat(\frac{x^{-a}}{1})=0$, and since
$\nat(\frac{x^{-a}}{1})=\frac{x^{-a}}{1}$ in $K[\Sigma]_{C}$ is an invertible element, it follows that $\nat(y)=0$.

It follows from Proposition \ref{vanishinghelper} that
$L^{\pnt}( \Sigma(-a))_a$ is exact, because $a\not\in  |\Ncc_{-\Sigma(-a)}|$. Hence
$H_\mm^i(K[\Sigma])_a\iso H^{i-1}(\Hom_\ZZ(\mathcal{C}_{\lpnt}(\Gamma_{\str_\Sigma(-a)}),K))\otimes_K K(-a)$.
\end{proof}

Note that $\str_\Sigma(C)$ is a partially ordered set with respect to inclusion.
We denote by $\Delta(\str_\Sigma(C))$ the order complex of $\str_\Sigma(C)\setminus\{C\}$
which is the simplicial complex where the faces are the chains of $\str_\Sigma(C)\setminus\{C\}$.
Let
$\Tilde H^{i}(\Delta(\str_\Sigma(C)))$ be the simplicial cohomology of the simplicial complex $\Delta(\str_\Sigma(C))$
with coefficients in $K$.

\begin{lemma}
\label{interpreter}
With the above notation we have for $i\in \ZZ$ that
$$
\Tilde H_{i}(\Gamma_{\str_\Sigma(C)})
\cong
\Tilde H_{i-\dim C}(\Delta(\str_\Sigma(C)))
\text{ and }
\Tilde H^{i}(\Gamma_{\str_\Sigma(C)})
\cong
\Tilde H^{i-\dim C}(\Delta(\str_\Sigma(C))).
$$
\end{lemma}
\begin{proof}
It is well-known that there exist a rational pointed fan $\Sigma'$ in $\RR^{d-\dim C}$
such that the face poset of $\Sigma'$ is isomorphic to $\str_\Sigma(C)$.
In fact, let
$\rho\colon \RR^d \to \RR^{d}/ \RR C\cong \RR^{d-\dim C}$ be the natural projection map
and consider for $D \in \str_\Sigma(C)$ the cones $\overline D=\rho(D)$.
(See \cite[Section 3.1]{FO93} for details.)
Now one checks that, up to a sign and a homological shift by $\dim C$,
the complex $\mathcal{C}_{\lpnt}(\Gamma_{\str_\Sigma(C)})$ is
the augmented oriented chain complex of the regular cell complex $\Gamma_{\Sigma'}$.
Hence we get
$
\Tilde H_{i}(\Gamma_{\str_\Sigma(C)})
\cong
\Tilde H_{i-\dim C}(\Gamma_{\Sigma'})
$
for $i \in \ZZ$.
Now $\Delta(\str_\Sigma(C))$ is just the barycentric subdivision of $\Sigma'$.
It follows from \cite[Proposition 4.7.8]{BJetal99} and \cite[Theorem 6.2.3]{BH}
that
$\Tilde H_{j}(\Gamma_{\Sigma'}) \cong\Tilde H_{j}(\Delta(\str_\Sigma(C)))$
for $j \in \ZZ$.
Analogously one shows the desired isomorphism for cohomology. This concludes
the proof.
\end{proof}

As a corollary we get a direct proof of the main result of \cite{BBR05} applied to toric face rings.
For Stanley--Reisner rings the next formula is due to Hochster.

\begin{corollary}
\label{hochster}
Let $\Sigma$ be a rational pointed fan  in $\RR^d$ and $a \in \ZZ^d$.
Then
\begin{eqnarray*}
H_\mm^i(K[\Sigma])
&\iso&
\bigoplus_{C\in \Sigma}
\bigoplus_{a \in -\relint(C)}
\Tilde H^{i-1}(\Gamma_{\str_\Sigma(C)}) \otimes_K K(-a)\\
&\iso&
\bigoplus_{C\in \Sigma}
\bigoplus_{a \in -\relint(C)}
\Tilde H^{i-\dim C-1}(\Delta(\str_\Sigma(C))) \otimes_K K(-a)
\end{eqnarray*}
as $\ZZ^d$-graded $K$-modules.
\end{corollary}
\begin{proof}
By Proposition \ref{vanishinghelper} we have
$H_\mm^i(K[\Sigma])_a=0$ for $a\not\in  |\Ncc_{- \Sigma}|$.
Thus we only have to consider elements in $|\Ncc_{- \Sigma}|$.
The assertion follows now from Theorem \ref{thmlcx}
and Lemma \ref{interpreter}.
\end{proof}

This decomposition implies some corollaries like a Cohen--Macaulay criterion for
toric face rings. For results in this direction see \cite[Section 5]{BBR05}.

%
%
%
\section{The canonical module}
\label{canonicalmodule}

In this section we  compute the $\ZZ^d$-graded dualizing complex for a toric face ring $K[\Sigma]$.
If $K[\Sigma]$ is Cohen--Macaulay we determine the $\ZZ^d$-graded canonical module of $K[\Sigma]$.
This results are generalizations to those ones known for Stanley--Reisner rings
(e.g. see \cite{BH} or \cite{GR84,GR842} for details) and affine monoid rings (e.g. see \cite{BH} or \cite{I} for details).

At first observe that for a cone $C$ in a rational pointed fan $\Sigma$ the
$K$-algebra $K[-C\cap \ZZ^d]$ is a $\ZZ^d$-graded subring of $K[\Sigma]_C$.
Let $k=\dim \Sigma$.
We set
$$
D^t=\bigoplus_{e_C\in \Gamma_\Sigma^{t-1}}K[-C\cap \ZZ^d] \text{ for }t=0,\ldots,k
$$
and we  consider the complex
$$
\CD
D^{\pnt}(\Sigma): 0@>>> D^0 @>\partial>> D^1 @>>>\cdots @>>> D^{k-1} @>\partial>> D^k @>>> 0,
\endCD
$$
where the differential is induced by the injection $D^{\pnt}(\Sigma)\to L^{\pnt}(\Sigma)$.
Observe that for $a\in \ZZ^d$ we have
$D^{\pnt}(\Sigma)_a \iso \Hom_\ZZ(\mathcal{C}_{\lpnt}(\Gamma_{\str_\Sigma(-a)})[-1],K)$.
We conclude that $$H_\mm^i(K[\Sigma])_a\iso H^i(D^{\pnt}(\Sigma))_a.$$
So the injection $D^{\pnt}(\Sigma)\to L^{\pnt}(\Sigma)$ of complexes is a quasi-isomorphism.

Let $(R,\mm)$ be a $\ZZ^d$-graded local ring, i.e.\ $R$ is a $\ZZ^d$-graded ring,
$\mm$ the unique $\ZZ^d$-graded maximal ideal of $R$,
and we assume that $\mm$ is also a maximal ideal.
By a dualizing complex in the category of $\ZZ^d$-graded $R$-modules
we mean a complex of $\ZZ^d$-graded injective $R$-modules
$$
\CD
I_{\lpnt}\colon
0@>>> I_{k}@>>> I_{k-1}@>>> \cdots @>>> I_{1}@>>> I_{0}@>>> 0,
\endCD
$$
such that $I_{\lpnt}$ has finitely generated homology modules and
the complex $\Hom_R(K,I_{\lpnt})$ is quasi-isomorphic to a
certain homological shift of the complex $0\to K\to 0$.
(For a comprehensive treatment we refer the reader to Chapter 15 in \cite{FO98}.)
For $\ZZ^d$-graded $R$-modules $M,N$ and $a \in \ZZ^d$
we denote by $\Hom_R(M,N)_a$ the set of degree preserving
$R$-homomorphisms
from $M$ to $N(a)$. Then
$\Hom_R(M,N)=\bigoplus_{a\in \ZZ^d} \Hom_R(M,N)_a$ is a $\ZZ^d$-graded $R$-module.
Note that if $M$ is finitely generated, then the underlying module of
$\Hom_R(M,N)$ is isomorphic to the usual (non-graded) set of $R$-module homomorphisms
from $M$ to $N$.
Further,
we let $L_{\lpnt}(\Sigma)$ be the $K$-dual complex $\Hom_K(L^{\pnt}(\Sigma),K)$ and
$D_{\lpnt}(\Sigma)$ be the complex $\Hom_K(D^{\pnt}(\Sigma),K)$.
Note that
$
D_t=\bigoplus_{e_C\in \Gamma_\Sigma^{t-1}}K[C\cap \ZZ^d]\text{ for }t=0,\ldots,k,
$
and that the differential  $\pi$ of $D_{\lpnt}(\Sigma)$
is the canonical projection multiplied
with the incidence function on $\Sigma$.

\begin{theorem}
Let $\Sigma$ be a rational pointed fan in $\RR^d$.
The complex $L_{\lpnt}(\Sigma)$ is a dualizing complex for $K[\Sigma]$
in the category of $\ZZ^d$-graded $K[\Sigma]$-modules.
The complex $D_{\lpnt}(\Sigma)$ is quasi-isomorphic to  $L_{\lpnt}(\Sigma)$.
\end{theorem}
\begin{proof}All $L^i$ are flat $K[\Sigma]$-modules,
so $\cdot\otimes_{K[\Sigma]}L^i$ is an exact functor. We have that $\Hom_K(\cdot,K)$ is an exact functor
and
it follows that
$$
\Hom_{K[\Sigma]}(\cdot,L_i)=\Hom_{K[\Sigma]}(\cdot,\Hom_K(L^{i},K))=\Hom_K(\cdot\otimes_{K[\Sigma]} L^i,K)
$$
is also exact.
We deduce that all $L_i$ are $\ZZ^d$-graded injective $K[\Sigma]$-modules. If $C\not= 0$, then
$K\otimes_{K[\Sigma]}K[\Sigma]_C=0$. Since $\Hom_{K[\Sigma]}(K,L_{\lpnt}(\Sigma))=
\Hom_K(K \otimes_{K[\Sigma]} L^{\pnt}(\Sigma),K)$
we deduce that the complex $\Hom_{K[\Sigma]}(K,L_{\lpnt}(\Sigma))$ is isomorphic to $ 0\to K\to 0$.
Thus $L_{\lpnt}(\Sigma)$ is a dualizing complex for $K[\Sigma]$.

Finally, the quasi-iso\-morphism between $D_{\lpnt}(\Sigma)$ and $L_{\lpnt}(\Sigma)$
follows from the quasi-iso\-morphism between  $D^{\pnt}(\Sigma)$ and $L^{\pnt}(\Sigma)$
since $\Hom_K(\cdot,K)$ is exact.
\end{proof}

$\omega_{K[\Sigma]}=H_k(L_{\lpnt}(\Sigma) )$ is called the {\em canonical module} of $K[\Sigma]$
where $k=\dim K[\Sigma]$.
(This generalizes the usual definition, see  \cite[Chapter 15]{FO98}.)

\begin{corollary}
\label{simplified_DC}
Let $\Sigma$ be a rational pointed fan in $\RR^d$ and let $K[\Sigma]$ be a Cohen--Macaulay ring.
Then the complex of $\ZZ^d$-graded $K[\Sigma]$-modules
$$
\CD
 0@>>> \omega_{K[\Sigma]}@>>>D_k @>\pi>> D_{k-1} @>>>\cdots @>>> D_{1} @>\pi>> D_0 @>>> 0
\endCD
$$
is exact.
\end{corollary}

We call a $\ZZ$-grading on a toric face ring $K[\Sigma]$ {\em admissible},
if with respect to this grading
$K[\Sigma]=\bigoplus_{i \in \NN}K[\Sigma]_i$ is a positively graded $K$-algebra
such that $K[\Sigma]_0=K$ and all components $K[\Sigma]_i$ are direct sums of
finitely many $\ZZ^d$-graded components.

\begin{corollary}
Let $\Sigma$ be a rational pointed fan in $\RR^d$ with $\dim \Sigma=k$ such that
$K[\Sigma]$ is Cohen--Macaulay and there exists an admissible $\ZZ$-grading on $K[\Sigma]$.
Let $C_1,\ldots,C_t\in\Sigma$ be the $k$-dimensional cones.
Choose $a_j\in \relint(C_j) \cap \ZZ^d$ for $j=1,\ldots,t$,
and let the degree of $x^{a_j}$ be $g_j$ with respect to the $\ZZ$-grading.
Set $h=\lcm(g_1,\ldots,g_t)$.
Then there exists an $\ZZ$-graded embedding
$$
\CD
\omega_{K[\Sigma]}(-h)@>>>K[\Sigma].
\endCD
$$
\end{corollary}
\begin{proof}
We consider the composition of homomorphisms
$$
\omega_{K[\Sigma]}(-h)
\longrightarrow\bigoplus_{j=1}^t K[C_j\cap\ZZ^d](-h)
\longrightarrow\bigoplus_{j=1}^t (x^b: b\in  \relint(C_j)\cap \ZZ^d)
\longrightarrow K[\Sigma],
$$
where $(x^b: b\in  \relint(C_j)\cap \ZZ^d)$ is the ideal in $K[C_j\cap\ZZ^d]$
generated by all elements in $\relint(C_j)\cap \ZZ^d$ for all $j$.
The middle  map
is given by (componentwise) multiplication with
$\bigl((x^{a_1})^{\frac{h}{g_1}},\ldots,(x^{a_t})^{\frac{h}{g_t}}\bigr)$.
The left  homomorphism is injective by Corollary \ref{simplified_DC}.
The middle  homomorphism is a direct sum of injective   homomorphisms, so it is injective. The last one is injective
since $\relint(C_i)\cap\relint(C_j)=\emptyset$ for $i\not=j$.
\end{proof}

\begin{remark}
The assumption that there exists an admissible grading on $K[\Sigma]$
is essential. Bruns and Gubeladze constructed in \cite[Example 2.7]{BRGU01}
a toric face ring $K[\Sigma]$ such that there is not  an admissible grading on $K[\Sigma]$.
To minimize $h$, one has to choose the interior elements in suitable way.
In general there is not much to say. In special cases there exists good bounds.
For example for Stanley--Reisner rings there always exists an admissible grading such
that one can choose $h$ as the dimension of the fan associated to the given simplicial complex.
(See \cite{GR84}.)
\end{remark}

%
%
%
\section{Gorenstein toric face rings}
\label{gorensteinsection}
The goal of this section is to characterize Gorenstein toric face rings.
For this we first introduce the following numbers.
For a cone $C$ in a rational pointed fan $\Sigma$
let $f_i(\str_\Sigma(C))$ be the number of $i$-dimensional cones in $\str_\Sigma(C)$,
and
$f_i(\Delta(\str_\Sigma(C)))$ be the number of $i$-dimensional faces in
the simplicial complex $\Delta(\str_\Sigma(C))$.
Using Lemma \ref{interpreter} we compute
\begin{eqnarray*}
(-1)^{\dim C}\sum_{i=0}^{\dim \Sigma} (-1)^{i} f_i(\str_\Sigma(C))
&=&
\sum_{i=-1}^{\dim \Sigma-\dim C} (-1)^{i} \dim_K \Tilde H^{i+\dim C}(\Gamma_{\str_\Sigma(C)})\\
&=&
\sum_{i=-1}^{\dim \Delta(\str_\Sigma(C))} (-1)^{i} \dim_K \Tilde H^{i}(\Delta(\str_\Sigma(C)))\\
&=&
\sum_{i=-1}^{\dim \Delta(\str_\Sigma(C))} (-1)^{i}f_i(\Delta(\str_\Sigma(C))).
\end{eqnarray*}
The latter number is nothing else that the {\em reduced Euler characteristic}
$\Tilde\chi_\Sigma (\Delta(\str_\Sigma(C)))$
of the simplicial complex $\Delta(\str_\Sigma(C))$ and
this motivates to define for
$C \in \Sigma$:
$$
\Tilde\chi_\Sigma (C)
=
(-1)^{\dim C}\sum_{i=0}^{\dim \Sigma} (-1)^{i} f_i(\str_\Sigma(C)).
$$
For a $\ZZ^d$-graded ring $R$ we let
$\Supp(R,\ZZ^d)=\{a \in \ZZ^d:R_a\not= 0\}$ be the {\em support} of $R$ in $\ZZ^d$.
\begin{theorem}
\label{gorensteinfirst}
Let $\Sigma$ be a rational pointed fan in $\RR^d$.
Then  $K[\Sigma]$ is Gorenstein if and only if $K[\Sigma]$ is Cohen--Macaulay and  there
exists $\sigma\in\big [\bigcap_{C \in \Sigma \text{ maximal}}C\big ]\cap\ZZ^d$ such that
we have for all $C \in \Sigma$
$$
\Tilde\chi_\Sigma (C)
=
\begin{cases}
(-1)^{\dim \Sigma -\dim C} , &\text{ if } \relint(C)\cap \ZZ^d -\sigma \subset \Supp(K[\Sigma],\ZZ^d),\\
0, &\text{else}.
\end{cases}
$$
In particular, if $\sigma=0$, then the latter condition reduces to
$\Tilde\chi_\Sigma (C)=(-1)^{\dim \Sigma -\dim C}$ for all $C \in \Sigma$. Note  that if  the above equivalent conditions are satisfied, then $\omega_{K[\Sigma]} \cong K[\Sigma](-\sigma)$.
\end{theorem}

We need some further notation.
For $a\in \ZZ^n$ we
let $f_i(a)$ be the number of $i$-dimensional cones in $\str_\Sigma(a)$.
Note that if $a\in  |\Ncc_{\Sigma}|$,
then $a \in \relint(C)$ for a unique $C \in \Sigma$.
Then $\str_\Sigma(a)=\str_\Sigma(C)$. Hence
$f_i(a)=f_i(\str_\Sigma(C))$ for all $i$ and we set $\Tilde\chi_\Sigma (a)=\Tilde\chi_\Sigma (C)$ in this case.
If $a\not\in  |\Ncc_{\Sigma}|$,
then $\str_\Sigma(a)=\emptyset$ and thus $f_i(a)=0$. We set $\Tilde\chi_\Sigma (a)=0$ in this case, and
we get $\Tilde\chi_\Sigma (a)=(-1)^{\dim C}  \sum_{i=0}^{\dim \Sigma} (-1)^{i} f_i(a)$ for all $a\in\ZZ^d$.

Let $M$ be a $\ZZ^d$-graded $K[\Sigma]$-module such that $\dim_K M_a <\infty$
for all $a \in \ZZ^d$.
We denote by $H_M(t)=\sum_{a \in \ZZ^n} \dim_K M_a t^a$ the {\em fine Hilbert series} of $M$.

\begin{proof}
Let $k=\dim \Sigma$ and assume that $K[\Sigma]$ is Gorenstein.
Since  $K[\Sigma]$ is Gorenstein, $K[\Sigma]$ is Cohen--Macaulay and there exists $\sigma\in\ZZ^d$ such that
$\omega_{K[\Sigma]} \cong K[\Sigma](-\sigma)$. From Corollary \ref{simplified_DC} we get an $\ZZ^d$-graded embedding
$$
\psi: K[\Sigma](-\sigma)\longrightarrow \bigoplus_{C \in \Sigma \text{ maximal}}K[C\cap \ZZ^d].
$$
Then $\psi(1)\in\bigoplus_{C \in \Sigma \text{ maximal}}(K[C\cap \ZZ^d])_\sigma$. Assume that there exists a maximal cone $C_0\in \Sigma$ such that
$K[C_0\cap \ZZ^d]_\sigma=0$, and take $a_0\in \relint(C_0)\cap\ZZ^d$. Then $0=x^{a_0}\psi(1)=\psi(x^{a_0})$, in contradiction with the fact that $\psi$ is injective.
Hence $\sigma$ belongs to all maximal cones in $\Sigma$.

Next  we have
${\omega_{K[\Sigma]}}_{a}\iso K[\Sigma]_{a-\sigma}$ for all $a\in \ZZ^d$.
Let $C \in \Sigma$ and $a \in \relint(C)\cap \ZZ^d$. Note that then $\str_\Sigma(a)=\str_\Sigma(C)$.
Now the exact complex in Corollary \ref{simplified_DC} in degree $a$ is
isomorphic to the exact complex of $K$-vector spaces
$$
0\to K[\Sigma]_{a-\sigma}\to K^{f_{k}(\str_\Sigma(C))} \to K^{f_{k-1}(\str_\Sigma(C))} \to\cdots \to K^{f_{1}(\str_\Sigma(C))} \to K^{f_{0}(\str_\Sigma(C))} \to 0.
$$
Since the alternating sum of the dimensions must be $0$, we get
$$
\dim_K K[\Sigma]_{a-\sigma}
=
(-1)^{k}\sum_{i=0}^{k} (-1)^if_{i}(\str_\Sigma(C))
=
(-1)^{\dim \Sigma -\dim C} \Tilde\chi_\Sigma (C).
$$
Observe that the right hand side does not depend on the chosen $a \in \relint(C)\cap \ZZ^d$.
The  conclusion on $\Tilde\chi_\Sigma (C)$ follows since $\dim_K K[\Sigma]_{a-\sigma} \leq 1$
depending whether $\relint(C)\cap \ZZ^d-\sigma$ is a subset of $\Supp(K[\Sigma],\ZZ^d)$ or not.
\smallskip

For the converse note that $K[\Sigma]$ Cohen--Macaulay implies the existence of a canonical module
$\omega_{K[\Sigma]}$ by Corollary \ref{simplified_DC}.
We show that $\omega_{K[\Sigma]} (\sigma)\cong K[\Sigma]$, which implies that $K[\Sigma]$ is Gorenstein.
First  we claim that the fine Hilbert series of $\omega_{K[\Sigma]}(\sigma)$ coincides with the one of
$K[\Sigma]$.
For $a \in \ZZ^d$,
the exact sequence of Corollary \ref{simplified_DC} in degree $a$ is isomorphic to
the sequence of $K$-vector spaces
$$
0\to {\omega_{K[\Sigma]}}_{a}\to K^{f_{k}(a)} \to K^{f_{k-1}(a)} \to\cdots \to K^{f_{1}(a)} \to K^{f_{0}(a)} \to 0.
$$
Since the alternating sum of the dimensions must be $0$, it follows
$$
\dim_K {\omega_{K[\Sigma]}}_{a}
=
(-1)^k \sum_{i=0}^k (-1)^i f_{i}(a)
=
\begin{cases}
(-1)^{\dim \Sigma -\dim C}\Tilde\chi_\Sigma (C) & \text{ if } a \in \relint(C), \text{ } C \in \Sigma,\\
0 & \text{else}.
\end{cases}
$$
If $a \not \in |\Ncc_{\Sigma}|$, then we saw that $\Tilde\chi_\Sigma (a)=0$.
Trivially $K[\Sigma]_a=0$ in this case,
because the set $|\Ncc_{\Sigma}|$ is the support of $K[\Sigma]$, which implies $K[\Sigma]_{a-\sigma}=0$, since $x^\sigma$ is not a zero--divisor.
If $a \in |\Ncc_{\Sigma}|$ we obtain $(\omega_{K[\Sigma]})_{a} \iso K[\Sigma]_{a-\sigma}$
as $K$-vector spaces. All in all we conclude that
the fine Hilbert-series $H_{\omega_{K[\Sigma]}}(t)=t^{\sigma}H_{K[\Sigma]}(t)$ coincide,
or equivalently, $H_{\omega_{K[\Sigma]}(\sigma)}(t)=H_{K[\Sigma]}(t)$.

Now we are ready to show $\omega_{K[\Sigma]} (\sigma)\cong K[\Sigma]$. There exists a non-zero $\ZZ^d$-graded homomorphism
$\phi:K[\Sigma]\to \omega_{K[\Sigma]}(\sigma)$ of degree $0$
induced by mapping $1$ to some $0\neq m^0 \in \omega_{K[\Sigma]}(\sigma)_0$.
We claim that $\phi$ is an isomorphism.
For this we consider the exact sequence
$$
\CD
0@>>>\Ker\phi @>>>K[\Sigma] @>\phi>>\omega_{K[\Sigma]}(\sigma)@>>>\Coker\phi @>>>0.
\endCD
$$
Since $K[\Sigma]$ and $\omega_{K[\Sigma]}(\sigma)$ have the same fine Hilbert series,
also the fine Hilbert series
of $\Ker\phi$ and $\Coker\phi$ coincide.
Let $m^0,m^1,\ldots,m^n$ be homogeneous generators of $\omega_{K[\Sigma]}(\sigma)$ such that
$\phi(1)=m^0$, and such that the images of $m^1,\ldots,m^n$ under projection
are a minimal system of generators of $\Coker\phi$.

We claim that $\omega_{K[\Sigma]}(\sigma)=K[\Sigma]m^0\oplus(K[\Sigma]m^1+\ldots+K[\Sigma]m^n)$. Observe that $\omega_{K[\Sigma]}$ is an indecomposable $K[\Sigma]$-module.
(For example use the fact that $\dim_K \Ext^i_{K[\Sigma]}(K,\omega_{K[\Sigma]})=\delta_{ik}$.)
Since $\omega_{K[\Sigma]}(\sigma)$ is indecomposable, it follows then $\omega_{K[\Sigma]}(\sigma) \cong K[\Sigma]$ because $m^0\neq 0$.

Since $\dim_K K[\Sigma]_a\le 1$ for all $a\in \ZZ^d$, we deduce
$$
\{a\in\ZZ^d:\Ker\phi_a\not=0\}\cup \{a\in\ZZ^d:K[\Sigma]_a m^0 \not=0\}=\Supp(K[\Sigma],\ZZ^d),
$$
$$
\{a\in\ZZ^d:\Ker\phi_a\not=0\}\cap \{a\in\ZZ^d:K[\Sigma]_a m^0\not=0\}=\emptyset.
$$
Now assume $K[\Sigma]m^0\cap (K[\Sigma]m^1+\ldots+K[\Sigma]m^n)\not=0$. Then there
exists a homogeneous element $0\neq y\in K[\Sigma]m^0\cap (K[\Sigma]m^1+\ldots+K[\Sigma]m^n)$.
It follows that there is a homogeneous element $x^b\in K[\Sigma]_b$ and some $m^i\in \omega_{K[\Sigma]}(\sigma)_c$, $i\ge 1$ such that
$b+c=\deg y$ and $x^bm^i\not=0$.
Since $\omega_{K[\Sigma]}(\sigma)_{\deg y}=\omega_{K[\Sigma]}(\sigma)_{b+c}\not=0$, we deduce $K[\Sigma]_{b+c}\not=0$.
On the other hand, $c\in\{a\in\ZZ^d:\Coker\phi_a\not=0\}= \{a\in\ZZ^d:\Ker\phi_a\not=0\}$
implies that $\Ker\phi_c=K[\Sigma]_c\not=0$.

From Corollary \ref{simplified_DC} we get a $\ZZ^d$-graded embedding
$$
\psi:\omega_{K[\Sigma]}\longrightarrow \bigoplus_{C \in \Sigma \text{ maximal}}K[C\cap \ZZ^d].
$$
If we assume $x^b x^c=0$ in $K[\Sigma]$, then $x^b x^c x^\sigma=0$.
Since $x^\sigma$ belongs to all maximal cones in $\Sigma$, we have $x^cx^\sigma=x^{c+\sigma}$. We deduce that
$x^b\big (\bigoplus_{C \in \Sigma \text{ maximal}}K[C\cap \ZZ^d]\big )_{c+\sigma}=0$. Since $\psi$ is injective it follows that
$x^b(\omega_{K[\Sigma]})_{c+\sigma}=x^b\omega_{K[\Sigma]}(\sigma)_c=0$,
in contradiction with $x^bm^i\not=0$. We deduce that $x^b x^c\not=0$, so
we have $\Ker\phi_{b+c}=K[\Sigma]_bK[\Sigma]_c=K[\Sigma]_{b+c}\not=0$.
Hence $\Im\phi_{b+c}=0$ which is a contradiction to the fact that $0 \neq y \in K[\Sigma]m^0=\Im\phi_{b+c}$.
This concludes the proof.

If in particular $\sigma=0$, then we always have $\relint(C)\cap \ZZ^d \subset \Supp(K[\Sigma],\ZZ^d)$.
\end{proof}

Let $\Sigma$ be a rational pointed fan in $\RR^d$
such that the conditions in Theorem \ref{gorensteinfirst} are satisfied. Then $K[\Sigma]$ is Gorenstein  and $\omega_{K[\Sigma]} \cong K[\Sigma](-\sigma)$.
A natural question is what one can say about $\sigma$. From the Theorem \ref{gorensteinfirst} we know that $\sigma$ belongs to all maximal cones in $\Sigma$.
Since $(\omega_{K[\Sigma]})_{\sigma}=K[\Sigma]_{\sigma-\sigma}=K[\Sigma]_{0}\not=0$,
we obtain from the exact sequence
$$
0\to (\omega_{K[\Sigma]})_{\sigma}\to K^{f_{k}(\str_\Sigma(\sigma))} \to K^{f_{k-1}(\str_\Sigma(\sigma))} \to\cdots \to K^{f_{1}(\str_\Sigma(\sigma))} \to K^{f_{0}(\str_\Sigma(\sigma))} \to 0.
$$
that, as in the proof above,
$\Tilde\chi_\Sigma (\sigma)=\pm 1$.

If $\str_\Sigma(\sigma)=\Sigma$ (i.e.\ $\sigma=0$),
then $\omega_{K[\Sigma]}\iso K[\Sigma]$ and we saw already that
in this case $\Tilde\chi_\Sigma (C)=(-1)^{\dim \Sigma - \dim C}$ for all $C \in \Sigma$.

Now consider $\str_\Sigma(\sigma)\not=\Sigma$ (i.e.\ $\sigma\not=0$).
Let $C\in \str_\Sigma(\sigma)$ and $a\in\relint(C)$.
Then also $a+\sigma\in\relint(C)$.
Since  $K[\Sigma]_{a+\sigma-\sigma}=K[\Sigma]_{a}\not=0$
we have $\Tilde\chi_\Sigma (a+\sigma)=(-1)^{\dim \Sigma -\dim C}$.
But $\Tilde\chi_\Sigma (\cdot)$ is constant on $\relint(C)$ and we deduce
$$
\Tilde\chi_\Sigma (C)=(-1)^{\dim \Sigma -\dim C}
\text{  and  }
\relint(C)\cap \ZZ^d -\sigma \subset \Supp(K[\Sigma],\ZZ^d)
\text{  for all  } C\in \str_\Sigma(\sigma).
$$

Let $C\in \Sigma(\sigma)$ and $a\in\relint(C)$.
Then $K[\Sigma]_a\not\subset \qq_{\Sigma(\sigma)}$,
while $K[\Sigma]_\sigma \subset \qq_{\Sigma(\sigma)}$.
If we suppose $K[\Sigma]_{a-\sigma}\not=0$,
then $K[\Sigma]_a=K[\Sigma]_{\sigma}K[\Sigma]_{a-\sigma}\subset \qq_{\Sigma(\sigma)}$,
which is a contradiction.
We deduce
$$
\Tilde\chi_\Sigma (C)=0 \text{  and  }
\relint(C)\cap \ZZ^d -\sigma \cap \Supp(K[\Sigma],\ZZ^d)=\emptyset
\text{  for all  }
C\in \Sigma(\sigma).
$$
All in all we get  the following equivalent reformulation of Theorem \ref{gorensteinfirst}.

\begin{theorem}
\label{gorensteinsecond}
Let $\Sigma$ be a rational pointed fan in $\RR^d$.
Then  $K[\Sigma]$ is Gorenstein if and only if $K[\Sigma]$ is Cohen--Macaulay and there
exists $\sigma\in\big [\bigcap_{C \in \Sigma \text{ maximal}}C\big ]\cap\ZZ^d$ such that:
\begin{enumerate}
\item
$$
\bigcup_{C\in \str_\Sigma(\sigma)}\relint(C)\cap \ZZ^d
=
\sigma+\Supp(K[\Sigma],\ZZ^d).$$
\item
For all cones $C \in \Sigma$
we have
$$
\Tilde\chi_\Sigma (C)=
\begin{cases}
(-1)^{\dim \Sigma - \dim C} &\text{ if } C \in \str_\Sigma(\sigma) ;\\
0 &\text{ if } C \in \Sigma(\sigma).
\end{cases}
$$
\end{enumerate}
\end{theorem}

\begin{remark}
Let $C \in \RR^d$ be a rational pointed cone.
For the affine monoid ring $K[C\cap \ZZ^d]$, which is also the toric face ring associated to
$\fan(C)$, Danilov and Stanley showed that the canonical module can be identified with the ideal generated
by $x^a$ for $a \in \relint(C)\cap \ZZ^d$. See \cite[Theorem 6.3.5]{BH} for details.
A corresponding result for Cohen--Macaulay Stanley--Reisner rings does not hold
with respect to its standard grading.
In fact, the canonical module may even not be identified with an $\ZZ^d$-graded
ideal of the given ring.
See \cite[Section 5.7]{BH} for a counterexample.

Now let $\Sigma$ be a rational pointed fan in $\RR^d$
and assume that $K[\Sigma]$ is Gorenstein. Then $\omega_{K[\Sigma]} \cong K[\Sigma](-\sigma)$.
Recall that $\Sigma(\sigma)=\Sigma \setminus \str_\Sigma(\sigma)$ is a subfan of $\Sigma$.
(Note that here $\sigma$ is not to be read as ``shift''.)
Consider the toric face ring $K[\Sigma(\sigma)]=K[\Sigma]/\qq_{\Sigma(\sigma)}$.
Observe that
$$
\Supp(\qq_{\Sigma(\sigma)}, \ZZ^d\}=\{a\in\bigcup_{C \in \str_\Sigma(\sigma)}\relint(C)\}.
$$
On the other hand
$x^\sigma$ is a non-zero divisor on $K[\Sigma]$ and
the ideal generated by $x^\sigma$ is
$(x^\sigma) \cong K[\Sigma](-\sigma)\cong \omega_{K[\Sigma]}$. It follows from \ref{gorensteinsecond} (i) that
$\Supp(\qq_{\Sigma(\sigma)}, \ZZ^d\}=\Supp((x^\sigma), \ZZ^d\}$,
so the two ideals coincide.
Hence
$\qq_{\Sigma(\sigma)}$ is a principle ideal generated by $x^\sigma$ and can be identified with
the canonical module $ \omega_{K[\Sigma]}$ of $K[\Sigma]$.
Note that $\qq_{\Sigma(\sigma)}=K[\Sigma]$ if $\sigma=0$.
\end{remark}

Similarly to the notion of an Euler simplicial complex we define:

\begin{definition}
Let $\Sigma$ be a rational pointed fan in $\RR^d$.
Then $\Sigma$ is called
an {\em Euler fan}
if $\Sigma$ is pure,
and $\Tilde\chi_\Sigma (C)=(-1)^{\dim \Sigma - \dim C}$
for all $C \in \Sigma$.
\end{definition}

\begin{corollary}
\label{needit}
Let $\Sigma$ be a rational pointed fan in $\RR^d$.
The following statements are equivalent:
\begin{enumerate}
\item
$K[\Sigma]$ is Cohen--Macaulay and $\Sigma$ is an Euler fan;
\item
$K[\Sigma]$ is
Gorenstein and $\omega_{K[\Sigma]} \cong K[\Sigma]$ as $\ZZ^d$-graded modules;
\item
$K[\Sigma]$ is Gorenstein and $\Tilde\chi_\Sigma (0)=(-1)^{\dim \Sigma}$.
\end{enumerate}
\end{corollary}
\begin{proof}
The equivalence of (i) and (ii) follows from Theorem \ref{gorensteinsecond}.

(ii) $\Rightarrow$ (iii):
Observe that $\{0\} \in \Sigma$. Then $\Tilde\chi_\Sigma (0)=(-1)^{\dim \Sigma}$ by
Theorem \ref{gorensteinfirst}.

(iii) $\Rightarrow$ (ii):
Let
$\omega_{K[\Sigma]} \cong K[\Sigma](-\sigma)$ as $\ZZ^d$-graded modules.
It follows from Theorem \ref{gorensteinsecond} that
$\Tilde\chi_\Sigma (0)=(-1)^{\dim \Sigma}$ if and only if $\{0\} \in \str_\Sigma(\sigma)$.
The latter is equivalent to the fact that $\str_\Sigma(\sigma)= \Sigma$ and thus $\sigma=0$.
\end{proof}

Now it remains the question what one can say for a
Gorenstein fan such that $\omega_{K[\Sigma]} \cong K[\Sigma](-\sigma)$ as $\ZZ^d$-graded modules for
a non-zero $\sigma$.

\begin{theorem}\label{gorenstein_noneuler}
Let $\Sigma$ be a rational pointed fan in $\RR^d$.
Then the following statements are equivalent:
\begin{enumerate}
\item
$K[\Sigma]$ is Gorenstein
and $\omega_{K[\Sigma]} \cong K[\Sigma](-\sigma)$ as $\ZZ^d$-graded modules
for some $0\not=\sigma \in \ZZ^d$;
\item
There exists an $0\not=\sigma\in\ZZ^d$ such that
$\Sigma(\sigma)$ is an Euler fan,
$x^\sigma$ is a non-zero divisor of $K[\Sigma]$ and
$K[\Sigma(\sigma)]=K[\Sigma]/(x^\sigma)$
is  Cohen--Macaulay.
\end{enumerate}
\end{theorem}
\begin{proof}
(i) $\Rightarrow$ (ii):
We know already that
$x^\sigma$ is a non-zero divisor of $K[\Sigma]$
and
$K[\Sigma(\sigma)]=K[\Sigma]/(x^\sigma)$.
Finally, $K[\Sigma(\sigma)]$ is Gorenstein by
\cite[Proposition 3.1.19]{BH} and
$\dim\Sigma(\sigma)=\dim\Sigma-1$.
Since $\sigma\neq 0$ we see by Theorem \ref{gorensteinsecond}
that $0=\Tilde\chi_{\Sigma}(0)=\Tilde\chi_{\Sigma(\sigma)}(0)\pm\Tilde\chi_{\Sigma}(\sigma)$,
which implies $\Tilde\chi_{\Sigma(\sigma)}(0) \neq 0$.
Thus $\Sigma(\sigma)$ is Euler by Corollary \ref{needit}.

(ii) $\Rightarrow$ (i):
It follows from \ref{needit} that $K[\Sigma(\sigma)]$ is Gorenstein.
Since $K[\Sigma(\sigma)]=K[\Sigma]/(x^\sigma)$ we obtain from the graded analogue of
\cite[Proposition 3.1.19]{BH} that $K[\Sigma]$ is Gorenstein.
\end{proof}

In addition to Theorem \ref{gorenstein_noneuler}, we make the following remark.
Let $\Sigma$ be a rational pointed fan in $\RR^d$.
Assume that $K[\Sigma]$ is Gorenstein
and $\omega_{K[\Sigma]} \cong K[\Sigma](-\sigma)$ as $\ZZ^d$-graded modules
for some $0\neq \sigma \in \ZZ^d$.
We claim that for $a\in\Supp(\qq_{\Sigma(\sigma)},\ZZ^d)$,
there are uniquely determined $b\in \Supp(K[\Sigma(\sigma)],\ZZ^d)$
and $n \in \NN$ such that $a=n\sigma+b$.

At first we prove the existence of an equation of the desired form.
Since $x^a \in \qq_{\Sigma(\sigma)}=(x^\sigma)$
there exists a maximal $n>0$ such that $x^a \in (x^{n\sigma})$ and
$x^a \notin (x^{(n+1)\sigma})$. For example one uses the fact that
otherwise $a-n\sigma \in \Supp(K[\Sigma], \ZZ^d)$ for all $n \in \NN$,
which is only possible if $-\sigma$ is an element of a maximal cone $C$ of $\Sigma$.
Since $\sigma \in C$ we obtain a contradiction because $C$ is pointed and $\sigma \neq 0$.
Thus $x^a=x^{n \sigma}x^b$ for some $b \in \Supp(K[\Sigma], \ZZ^d)$. By the choice of
$n$ we have $x^b \not\in \qq_{\Sigma(\sigma)}$ and therefore the residue class of
$x^b$ is not zero in $K[\Sigma(\sigma)]$ which means that $b \in \Supp(K[\Sigma(\sigma)],\ZZ^d)$.
Hence $a=n\sigma+b$ is of the desired form.

Let $a=m\sigma+c$ for some $m \in \NN$ and $c\in\Supp(K[\Sigma(\sigma)],\ZZ^d)$ and assume that $m\geq n$.
Then $x^b=x^{(m-n)\sigma}x^c$ is not zero in $K[\Sigma]$.
The residue class in $K[\Sigma(\sigma)]$ is also not zero, which is only possible if $n=m$ and $b=c$.
This implies the uniqueness of such an equation.

\begin{proposition}Let $\Sigma$ be a rational pointed fan in $\RR^d$. If $K[\Sigma]$ is Gorenstein then
either $\Sigma$ is an Euler fan with $K[\Sigma]$ Cohen--Macaulay, or there is an Euler fan $\Sigma'$ with $K[\Sigma']$ Cohen--Macaulay
such  that $K[\Sigma]$ is isomorphic as a $\ZZ^d$-graded $K$ vector space to
a polynomial ring in one variable $K[\Sigma'][z]$ with coefficients in
$K[\Sigma']$ and where $\deg z=\sigma$.
\end{proposition}

Observe that usually this is not an isomorphism of $K$-algebras.
E.g.\ let $\Sigma=\fan(C)$ be the fan associated to one rational pointed cone $C$ in $\RR^d$.
Then $\sigma \neq 0$ is an element of $\relint(C)$. Now
$K[\Sigma]$ is an integral domain, while $K[\Sigma(\sigma)]$
has zero divisors and can not be a subalgebra of $K[\Sigma]$, which would be the case
if the isomorphism above is a ring homomorphism in an obvious way.

\end{document}